\documentclass[11pt]{article}

\usepackage{graphicx} 
\usepackage{geometry}
\usepackage{graphicx}
\usepackage{dcolumn}
\usepackage{bm}

\usepackage[utf8]{inputenc}
\usepackage[T1]{fontenc}
\usepackage{mathptmx}
\usepackage{etoolbox,authblk}

\usepackage{graphicx,svg,gensymb,float}
\usepackage{newtxtext}
\usepackage{newtxmath}
\usepackage{wrapfig}
\usepackage{blindtext}

\usepackage{graphicx,wrapfig,tikz}
\usepackage{subcaption}
\usepackage[export]{adjustbox}
\usepackage{natbib}
\usepackage{hyperref,ulem}

\makeatletter

\newcommand{\framedgraphic}[1]{
        \centering
        \begin{tikzpicture}
            \node[inner sep=0pt] (image) at (0,0) {\includegraphics[width=\textwidth]{#1}};
            \draw[black, line width=1pt, opacity=0.1] (image.south west) rectangle (image.north east);
        \end{tikzpicture}
        \caption{}
}

\makeatletter
\def\@email#1#2{%
 \endgroup
 \patchcmd{\titleblock@produce}
  {\frontmatter@RRAPformat}
  {\frontmatter@RRAPformat{\produce@RRAP{*#1\href{mailto:#2}{#2}}}\frontmatter@RRAPformat}
  {}{}
}%
 
\title{\bf Control of the Fluidic Pinball using the Quadratic-Quadratic Regulator}

\date{}

\begin{document}

\author[1]{Ali Bouland}
\author[2]{Jeff Borggaard}
\affil[1]{Department of Mechanical Engineering, Virginia Tech, Blacksburg, VA}
\affil[2]{Department of Mathematics, Virginia Tech, Blacksburg, VA}

\footnotetext[1]{Corresponding author: {\tt bouland@vt.edu}.}

\maketitle

\begin{abstract}
The fluidic pinball presents a significant benchmark for nonlinear flow control, managing the complex interactions of three cylinder wakes. This study addresses the stabilization of the fluidic pinball to its unstable steady-state solution using a model-based nonlinear feedback strategy. We propose a framework that combines interpolatory model order reduction (IMOR) with the quadratic-quadratic regulator (QQR), a feedback control methodology that is specifically suited to the quadratic nonlinearity of the Navier-Stokes equations. A finite element model (FEM) of the problem coupled with IMOR is used to produce a reduced-order model (ROM) that accurately represents the input-output dynamics of the actuated wake.  The performance of the QQR control is evaluated against the traditional linear feedback control for two different Reynolds numbers, $Re_D = 30$ and $Re_D = 50$. At $Re_D = 30$, the QQR controller is able to stabilize the wake and reaches the desired performance criteria 40.1\% faster than using a linear feedback controller. More significantly, at $Re_D = 50$, the QQR controller successfully stabilizes the wake, whereas the linear controller fails to overcome the nonlinearity of the flow.  The QQR control effectively suppresses vortex shedding, resulting in the elimination of lift oscillations and a reduction in the drag coefficient.  These results demonstrate that the IMOR-QQR framework provides an effective model-based control strategy that can manage nonlinear hydrodynamic instabilities in such complex wake flows. 
\end{abstract}
\newcommand{\keywords}[1]{%
  \par\vspace{0.5\baselineskip}%
  \noindent{\textbf{Keywords:} #1}
}


\maketitle

\section{Introduction and Background}
\label{sec:headings}

\subsection{Flow Control Problems}
The problem of flow control has a rich history with applications that include turbulence control, increased lift and decreased drag, minimizing vortex-induced vibrations, and enhancing or preventing mixing. We consider fluids that are described by the Navier-Stokes equations (NSE), non-dimensionalized as
$$ \begin{aligned}
\frac{\partial \boldsymbol{v}}{\partial t}+\boldsymbol{v} \cdot \nabla \boldsymbol{v} & =-\nabla p + \nabla \cdot \left(\frac{1}{Re} \left(\nabla \boldsymbol{v}+\nabla \boldsymbol{v}^{\mathrm{T}}\right)\right), \\
 \nabla \cdot \boldsymbol{v} & =0,
\end{aligned}$$
where $\boldsymbol{v} = (v_x,v_y)$ is the velocity field, $p$ is the pressure field, $t$ is time, and $Re$ is the Reynolds number. Semi-discrete approximations to this system of partial differential equations lead to a large system of differential-algebraic equations of the form
\begin{equation} \label{eq:CFD_DAE} 
    \left[\begin{array}{cc}
        {E}_{11} & 0 \\
        0 & 0
        \end{array}\right]\left[\begin{array}{c}
        \dot{\boldsymbol{x}}_{1} \\
        \dot{\boldsymbol{x}}_{2}
        \end{array}\right]=\left[\begin{array}{cc}
        {A}_{11} & {A}_{21}^{T} \\
        {A}_{21} & 0
        \end{array}\right]\left[\begin{array}{c}
        \boldsymbol{x}_{1} \\
        \boldsymbol{x}_{2}
        \end{array}\right]+\left[\begin{array}{c}
        {B}_{1} \\
        {B}_{2}
        \end{array}\right] \boldsymbol{u} + 
        \left[\begin{array}{c}{N} (\boldsymbol{x}_1 \otimes \boldsymbol{x}_1) \\ 0 \end{array}\right] ,
\end{equation}
where $\boldsymbol{x}_1$ are velocity coefficients, $\boldsymbol{x}_2$ are pressure coefficients and the discretized operators are: ${E}_{11}$ (Gram matrix of the velocity basis), ${A}_{11}\boldsymbol{x}_1$ (viscous stress), ${A}_{21}^T\boldsymbol{x}_2$ (pressure gradient), ${A}_{21}\boldsymbol{x}_1$ (velocity divergence), ${N}(\boldsymbol{x}_1\otimes\boldsymbol{x}_1)$ (convective velocity term), and $\otimes$ is the Kronecker product.  The terms ${B}_1$ and ${B}_2$ are used to model flow actuation with control input $\boldsymbol{u}$.

For even modest values of the Reynolds number, the dimensions of $\boldsymbol{x}_1$ ($n_1$) and $\boldsymbol{x}_2$ ($n_2$) can be large.  This creates several difficulties when developing approximate solutions to flow control problems, due to the need to perform multiple simulations or solve nonlinear matrix equations with dimensions on the order of the system size.  These difficulties also include the complexity of simulating the controlled full-order system using computational fluid dynamics (CFD), as the stiffness introduced by the control terms can significantly increase the computational cost.  Thus, model order reduction (MOR) is usually used throughout this process.  Note that since the NSE is nonlinear, linear controllers may not be able to control the flow in many cases, especially when trying to stabilize a solution from flows that are far from those used to perform the linearization. This motivates the search for nonlinear control algorithms that can be applied to flow control problems.  In this paper, we investigate nonlinear feedback control strategies using the fluidic pinball benchmark problem \cite{cornejo2019artificial} with the goal of using cylinder rotation to stabilize the unstable steady-state solution.  We assess the performance of control strategies by simulating the full-order controlled system.

\subsection{History of the Pinball Problem}
 \subsubsection {One and Two Cylinders}

For a single cylinder with rotational oscillation, initial studies focused on the ability of rotational oscillations to influence the wake structure~\cite{okajima1975viscous,taneda1978visual}.  Subsequent studies considered the generated vortex patterns produced by this actuation; see \cite{williamson1988vortex}, \cite{fujisawa1998vortex}, and \cite{lee2006flow}.  Computational investigations, such as \cite{ou1991control,shiels1996computation}, confirmed its potential for drag reduction. In \cite{baek1998numerical}, a numerical investigation of vortex formation was performed while varying the oscillation frequency and the maximum angular amplitude. Several active control techniques for vortex suppression have been investigated in the literature, and a comprehensive summary can be found in \cite{doreti2018control}. A range of control actuation methods have been considered, including surface suction  \cite{chen2014experimental}, mass-flow \cite{rabault2019artificial}, synthetic jets \cite{wang2016control}, plasma actuators \cite{jukes2009flow}, and rotational oscillation of the cylinder \cite{homescu2002suppression}. The latter actuation mechanism is used in the pinball problem and is the emphasis of this paper. Since cylinder rotation does not inject flow into the domain, this simplifies \eqref{eq:CFD_DAE} since ${B}_2 = {0}.$

Other authors have investigated twin (side-by-side) cylinder flows: for example \cite{carini2014first,carini2014origin} investigate the instability and sensitivity of the flow past two side-by-side cylinders and the asymmetric unsteady wake that occurs behind them. In \cite{sarvghad2011numerical}, numerical simulations were performed to investigate the effect of different spacing distances and Reynolds numbers. On the other hand, the authors in \cite{kumar2011flow} perform an experimental study with rotating cylinders, while varying the Reynolds number, spacing distance, and rotation speed. In \cite{chen2022numerical}, the authors performed a numerical simulation with passive-suction-jet control to manipulate the vortex shedding. Others, such as \cite{chan2011vortex}, have investigated the control of vortex shedding by counter rotation of the two cylinders.

 \subsubsection {The Pinball Problem}
 
 The so-called pinball problem that we consider closely follows that of \cite{deng2018route}. The flow is assumed to be uniform from left to right, with speed $V_\infty$. The three cylinders form an equilateral triangle pointing upstream and all have the same diameter $D$, with a triangle radius of $1.5D$, measured from the centroids of the trio of cylinders. The flow is assumed to be two-dimensional with kinematic viscosity $\nu$. All quantities are non-dimensionalized with the velocity magnitude $V_\infty$, cylinder diameters $D$, and kinematic viscosity $\nu$ leading to the Reynolds number defined as $Re_D = { V_\infty D}/{\nu}$. The computational domain is defined as
 \begin{align*}
     \Omega = \{ &(x,y) \in \mathbb{R}^2: (x,y) \in [0,26] \times[0,12] \cap\\
     &(x-x_i)^2 + (y-y_i)^2 \geq 1/4,\ i=1,2,3 \}.
 \end{align*}
 The centers of the cylinders, ordered counter-clockwise, are located at 
 \begin{align*}
     x_1 &= 6-3/2 \operatorname{cos}30\degree, && y_1 = 6,\\
     x_2 &= 6, && y_2 = 6-3/4,\\
     x_3 &= 6, && y_3 = 6+3/4.
 \end{align*}
 Throughout the paper, the upstream, bottom, and top cylinders are referred to as the first, second, and third cylinders, respectively. 
 
 Some of the most notable early papers on the fluidic pinball problem include \cite{deng2018route}, which focuses on the dynamics and route to chaos of the unforced pinball problem. Another important paper, \cite{pastur2019reduced} explores potential reduced-order modeling approaches. In that paper, the authors derive a proper orthogonal decomposition (POD) reduced-order model, for Reynolds number $Re_D = 30$ and stationary cylinders. Another paper on reduced-order modeling is \cite{deng2020low}, where the authors propose a least order Galerkin model applied to the fluidic pinball. A relevant reference, \cite{cornejo2019artificial}, discusses the application of machine learning methods to control the fluidic pinball, with the goal of reducing net drag. Note that this is a model-free method and the control law is derived by solving a non-smooth optimization problem.  In \cite{maceda2021stabilization}, the authors explore machine learning control methods in three different search spaces, open loop, closed loop, and output feedback control. Some more recent examples include \cite{Deng2021Galerkin}, \cite{Li2022explorative}, \cite{haodong2023how},  and \cite{Marra2024self}. 
 
 \subsection{Wake Stabilization}

In this work, we consider using a model-based control approach that combines
interpolatory model reduction and the nonlinear quadratic--quadratic regulator
to drive the pinball flow to its steady-state solution. Stabilizing the unstable
steady-state solution also stabilizes the wake and exhibits desirable drag reduction.
We provide a high-level overview of the methodology below. There are several design
choices that we made for this study, and the following sections will motivate these
decisions.

\begin{figure}[ht]
    \centering
    \includegraphics[width=0.8\columnwidth]{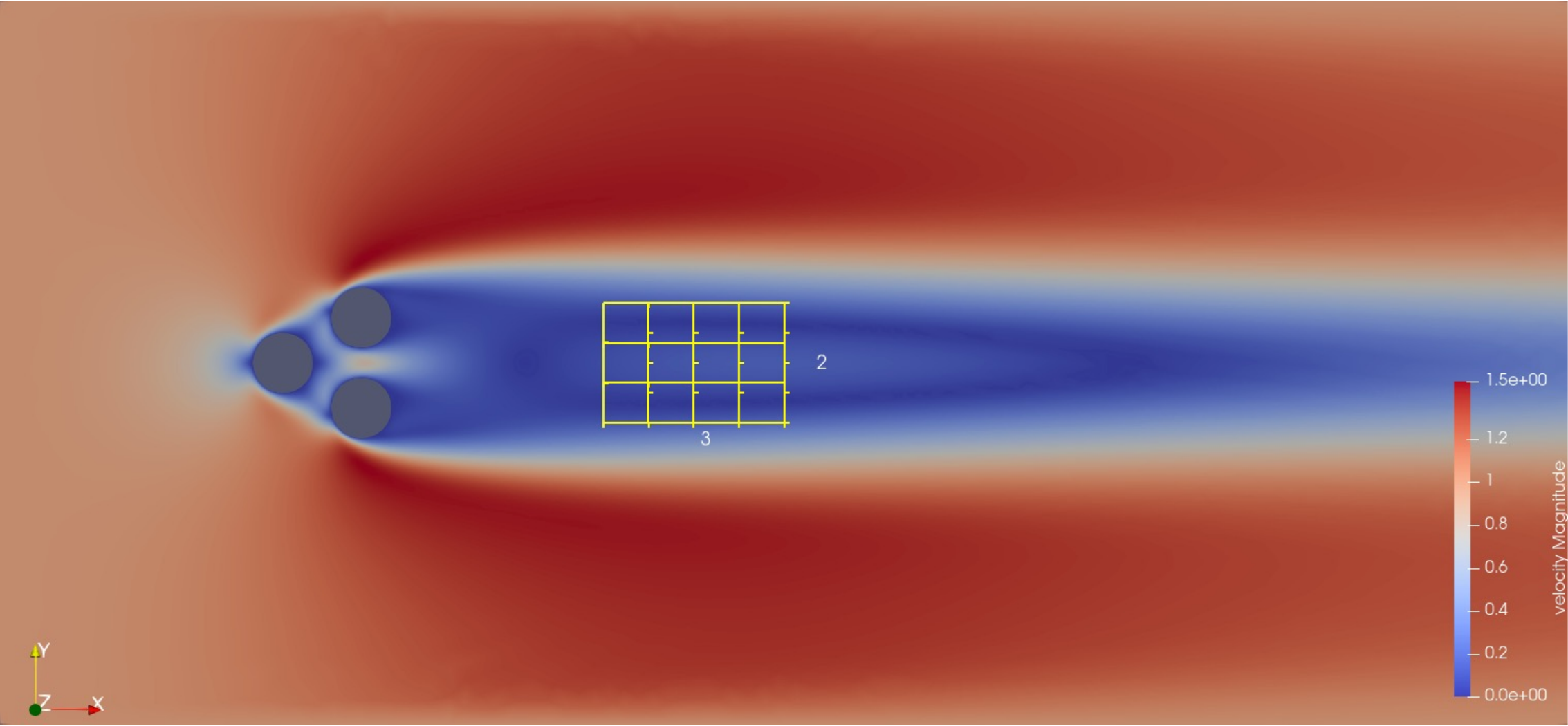}
    \caption{The velocity magnitude of the steady-state solution $\boldsymbol{v}_{ss}$. Boxes in the wake region indicate the regions where average velocity fluctuations are to be controlled.}
    \label{fig:Cmatrix}
\end{figure}

Let $(\boldsymbol{v}_{ss}, p_{ss})$ represent the steady-state solution to the pinball flow at a
given Reynolds number (see figure~\ref{fig:VorticityPlot} for the $Re_D=50$ case). Writing
\[
\boldsymbol{v} = \boldsymbol{v}_{ss} + \boldsymbol{v}', \qquad p = p_{ss} + p',
\]
and substituting into the Navier--Stokes equations leads to the following nonlinear
equations for the perturbation variables:
\begin{align}
\label{eq:fluctuatingNS}
\frac{\partial {\boldsymbol{v}'}}{\partial t}
+ \boldsymbol{v}_{ss} \cdot \nabla {\boldsymbol{v}'}
+ {\boldsymbol{v}'} \cdot \nabla \boldsymbol{v}_{ss}
+ {\boldsymbol{v}'} \cdot \nabla{\boldsymbol{v}'}
&= -\nabla p' +
\frac{1}{Re} \nabla \cdot \left( \nabla {\boldsymbol{v}'} 
+ \nabla {\boldsymbol{v}'}^{\mathrm{T}} \right),\\
 0 &= \nabla \cdot {\boldsymbol{v}'}
\end{align}

In Section~\ref{sec:CFDSpecific}, we detail a finite element approximation to these equations with discrete variables $\boldsymbol{x}_1$ and $\boldsymbol{x}_2$ to approximate $\boldsymbol{v}'$ and $p'$, respectively.
These equations satisfy a differential-algebraic equation (DAE) of the form of equation \eqref{eq:CFD_DAE}. 
Thus our model has exactly the form of the discrete Navier--Stokes equations, where
the ${A}_{11}$ matrix is updated to include the linear terms in (\ref{eq:fluctuatingNS}) that involve $\boldsymbol{v}_{ss}$.
We introduce sensors to measure spatial averages of the velocity
perturbation components in 12 regions in the wake for 24 controlled outputs.  These are located in a rectangle from $(10,5)$ to $(13,7)$ that includes the vortex formation length, cf.~\cite{wang2024vortex}, with the motivation that controlling the perturbations in this region effectively stabilizes the entire wake.  These measurements can be represented in the finite element model as
\begin{equation}
\boldsymbol{y}(t) = {C} \boldsymbol{x}_1(t).
\label{eq:OutPutEquation}
\end{equation}
The control problem now is to design the control inputs
$\boldsymbol{u}(t)\in \mathcal{L}_2(0,\infty; \mathbb{R}^3)$ representing the tangential velocities of the three cylinders 
that solve
\begin{equation}
\min_{\boldsymbol{u}} J(\boldsymbol{u})
\;=\;
\int_0^\infty \left( \|\boldsymbol{y}(t)\|^2 + \boldsymbol{u}(t)^T  {R}\,\boldsymbol{u}(t) \right)\,dt,
\label{eq:minJ}
\end{equation}
where $ {R}>0$ is a symmetric weighting matrix and $\boldsymbol{y}$ depends on $\boldsymbol{u}$ through the
solution of \eqref{eq:CFD_DAE} and \eqref{eq:OutPutEquation}.

As noted in the introduction, reduced-order models are an essential element in
model-based fluid control. The above setting has been designed to have
low-dimensional inputs and outputs and enables model reduction. The methodology we use is motivated in Section~\ref{sec:MOR}.  
Model reduction methods for fluids typically use reduced basis functions
\[
\mathbf{\Phi} = [\boldsymbol{\phi}_1,\dots,\boldsymbol{\phi}_r],
\]
satisfying ${{A}}_{21}\boldsymbol{\phi}_i=0$ for each $i=1,\ldots, r$. Thus, we write
\begin{equation}\label{eq:reducedApprox}
\boldsymbol{x}_1 \approx \boldsymbol{\Phi} \boldsymbol{x}_r.
\end{equation}
and perform Galerkin or Petrov--Galerkin projection with \eqref{eq:CFD_DAE} to find a controlled dynamical system
for $\boldsymbol{x}_r$ of the form
\begin{align}
\boldsymbol{\dot{x}}_r &= {A}_r \boldsymbol{x}_r + {B}_r \boldsymbol{u} + {N}_r(\boldsymbol{x}_r \otimes \boldsymbol{x}_r), \label{eq:RomState}\\
\boldsymbol{y}_r &= {C}_r \boldsymbol{x}_r.  \label{eq:RomOutPut}
\end{align}
In this work, we use IRKA on the linear portion of \eqref{eq:CFD_DAE} to develop $\mathbf{\Phi}$. This strategy incorporates the effects of the
inputs and outputs of the system. 

The reduced control problem is then to minimize ${J}_r(\boldsymbol{u})$ where ${J}_r$ matches \eqref{eq:minJ} except the term $\boldsymbol{y}(t)$ is replaced with $\boldsymbol{y}_r(t)$. The output $\boldsymbol{y}_r(t)$ depends on $\boldsymbol{u}$ through the solution to \eqref{eq:RomState} and \eqref{eq:OutPutEquation}. Minimizing the quadratic cost $J_r$ subject to the quadratic state equations \eqref{eq:RomState} is a quadratic--quadratic regulator problem. The optimal control can be represented as a solution to the Hamilton-Jacobi-Bellman (HJB) equation, which is a nonlinear equation in $r$ spatial dimensions. 
As we show in section~\ref{sec:polynFeedback}, the optimal control can be written in feedback form ${\boldsymbol{u}} = \mathcal{K}(\boldsymbol{x}_r)$
and approximated using a polynomial in Kronecker-product form, for example
\[
\boldsymbol{u}(\boldsymbol{x}_r) = \boldsymbol{k}_1 \boldsymbol{x}_r + \boldsymbol{k}_2 (\boldsymbol{x}_r \otimes \boldsymbol{x}_r).
\]

While our original objective is to minimize the full-order cost subject to the
full-order dynamics \eqref{eq:CFD_DAE}, the control low is developed for the reduced model.  Thus in section~\ref{sec:closedLoop} we lift the controller back
to the full-order states. The control performance is evaluated using the full-order model measured using $J(\boldsymbol{u})$.  The discussion in the following details our methodology.

\section{Detailed Problem Description}
\label{sec:ProblemDescription}

In the following sections, we briefly survey the important elements of the wake stabilization problem and motivate our design choices.  These include using FEniCS to perform CFD simulations, using interpolatory model reduction (IMOR) methods instead of more common methods for reduced-order modeling of fluids, incorporating quadratic terms in our feedback control design using QQR, and integrating these strategies.

    \subsection{CFD Specifics}
    \label{sec:CFDSpecific}
        \subsubsection{FEniCS and details}
\begin{figure} 
            \centering
            \includegraphics[width=0.8\columnwidth]{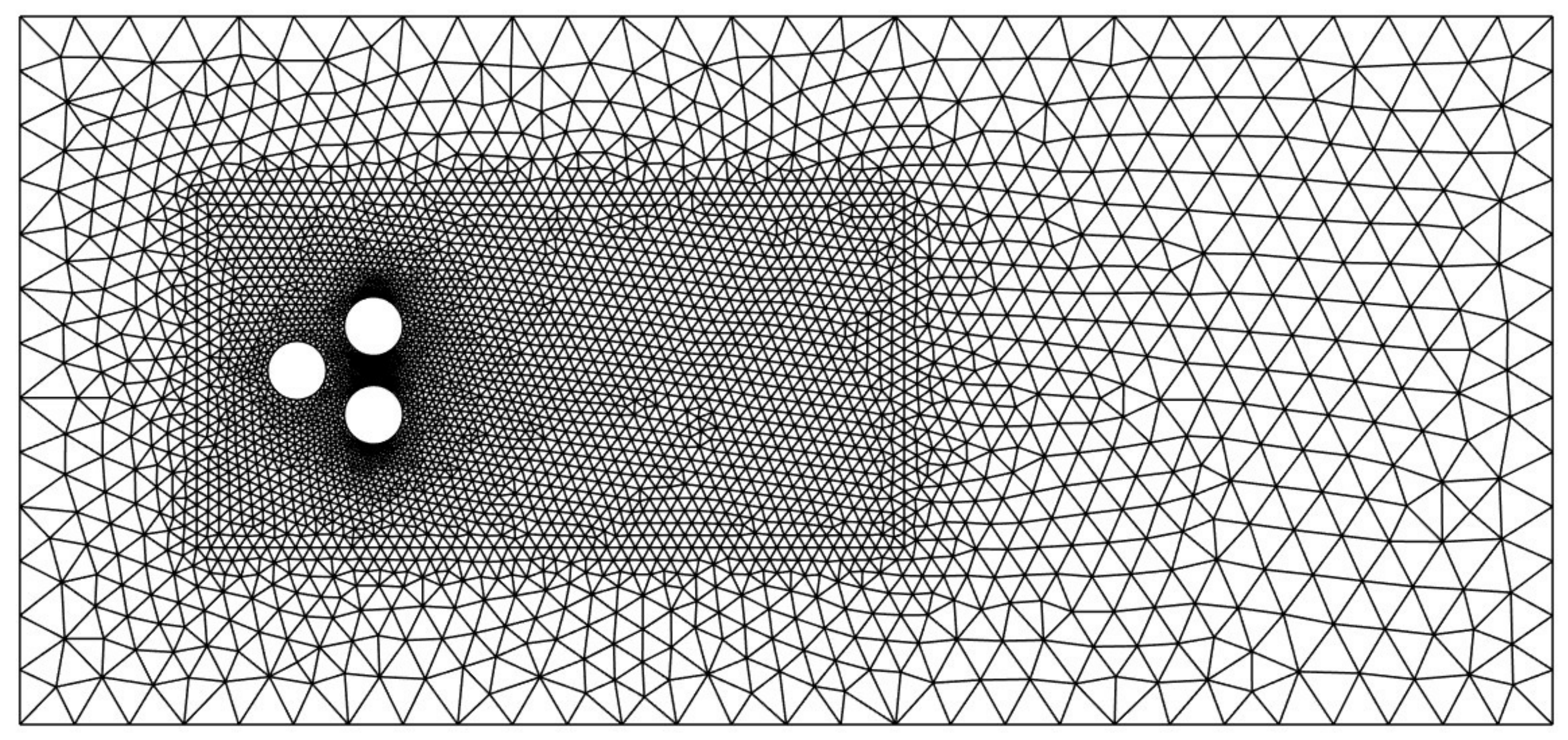} 
            \caption{Plot of the mesh used for the fluidic pinball.} \label{fig:MeshFigure}
        \end{figure}

        The finite element simulation was carried out using the FEniCS library, see \cite{logg2010dolfin}. The mesh was generated using the open source package Gmsh  \cite{geuzaine2009gmsh}. The pinball domain $\Omega$ is meshed using $7008$ P2-P1 elements with $14199$ nodes clustered near the cylinders and wake region (see figure~\ref{fig:MeshFigure}).  A constant timestep of $\Delta t=0.01$ is used for time-dependent flow simulations. 
        
         To find the steady--state velocity field $\boldsymbol{v}_{ss}$ required to build the ${A}_{21}$ matrix, the steady Navier--Stokes equations \eqref{eq:CFD_DAE} are solved using a Newton-based nonlinear solver. The open-loop (no control) and closed-loop simulations use the Incremental Pressure Correction Scheme (IPCS). The IPCS scheme consists of three main steps \cite{goda1979multistep,langtangen2017solving}. 

         First of all, a prediction of the velocity at the next step $n+1$, denoted by $\boldsymbol{v}^*$ is found using the weak form
         
     \[\begin{aligned}  \left ( \frac{\boldsymbol{v}^{*}-\boldsymbol{v}^n}{\Delta t}, \mathbf{w}\right )+\left ( \boldsymbol{v}^n \cdot \nabla \boldsymbol{v}^n, \mathbf{w}\right )+\left (\sigma\left(\boldsymbol{v}^{n+\frac{1}{2}}, p^n\right), \frac{1}{2}(\nabla \mathbf{w} +  \nabla \mathbf{w}^T)\right ) \\  +\left ( p^n  \hat{\boldsymbol{n}}, \mathbf{w}\right )_{\partial \Omega}-\left (\mu \nabla \boldsymbol{v}^{n+\frac{1}{2}} \cdot \hat{\boldsymbol{n}}, \mathbf{w}\right )_{\partial \Omega}&= \mathbf 0,\end{aligned}\]
where $\boldsymbol{v}^{n+\frac{1}{2}} = \frac{\boldsymbol{v}^n + \boldsymbol{v}^*}{2}$,  $\sigma(\boldsymbol{v},p) = \mu (\nabla \boldsymbol{v} + \nabla \boldsymbol{v}^T) - p\,\mathbf{I}$, $\mathbf{w}$ is a test function, and the inner products $(a,b) \equiv \int_{\Omega} ab \operatorname{d\Omega} $ and $(a,b)_{\partial \Omega} \equiv \int_{\partial \Omega} ab \operatorname{ds}$. To impose the divergence free condition, we advance the pressure using the pressure Poisson equation for $p^n$, then correct the velocity by subtracting the gradient of the pressure increment.  The weak form of the pressure update step is the following:
    \[\left (\nabla p^{n+1}, \nabla q\right )=\left (\nabla p^n, \nabla q\right )-\frac{1}{\Delta t}\left (\nabla \cdot \boldsymbol{v}^*, q\right ),\]
     where $q$ is a test function from the pressure space. Finally, we compute the corrected $\boldsymbol{v}^{n+1}$ by solving
    \[\left (\left(\boldsymbol{v}^{n+1}-\boldsymbol{v}^*\right), \mathbf{w}\right )=-\Delta t\left (\nabla\left(p^{n+1}-p^n\right), \mathbf{w}\right ).\]
     
        \begin{figure}[!htb]
            \centering
            \includegraphics[width=0.8\columnwidth]{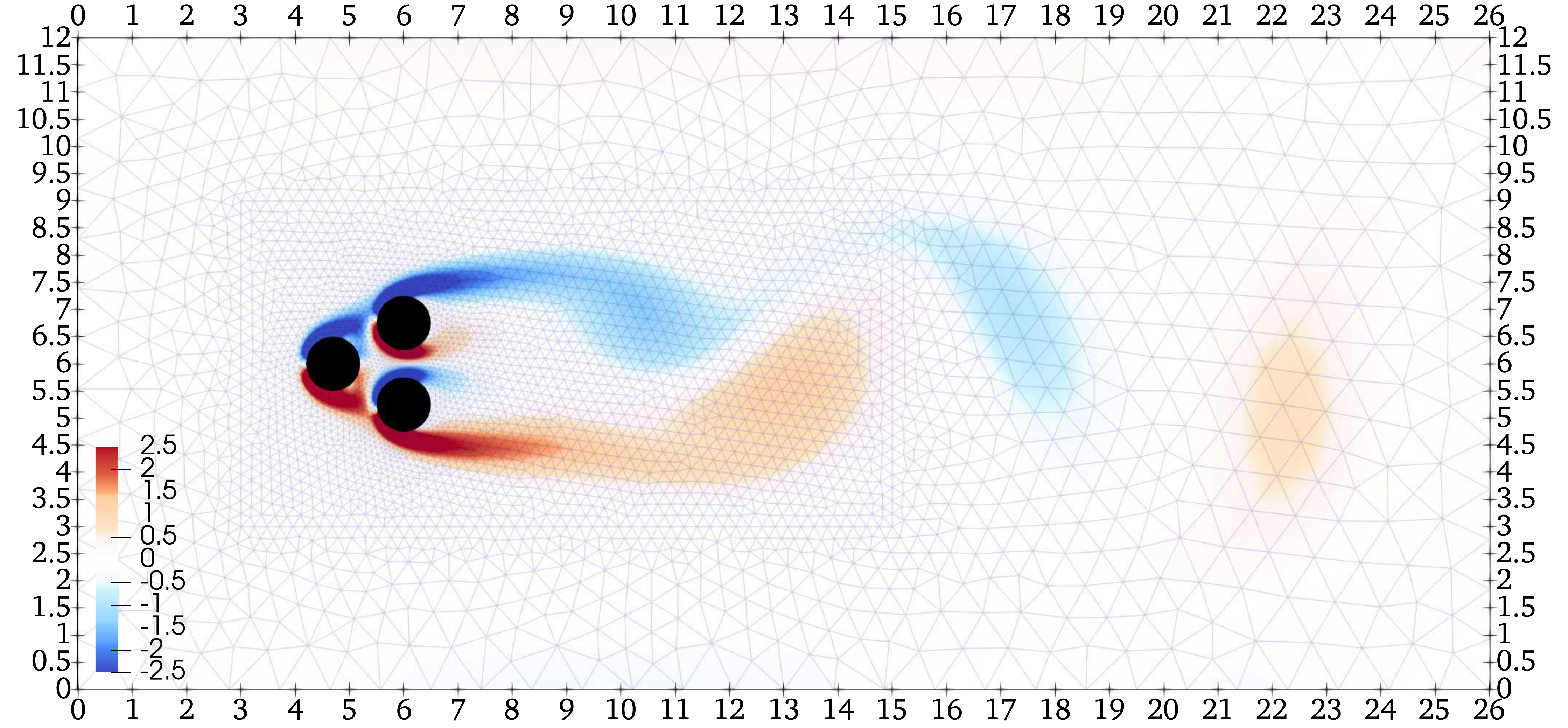} 
            \caption{Vorticity Plot for $Re_D=50$.} \label{fig:VorticityPlot}
        \end{figure}

    \subsubsection{Boundary Conditions}
    The boundary conditions closely follow what was used in \cite{deng2018route}. They consist of Dirichlet boundary conditions for the three cylinders ($\boldsymbol{v} = 0$), free stream boundary conditions for the inlets and walls ($\boldsymbol{v} = \boldsymbol{e}_x$). For the outlet, stress-free boundary conditions are assumed. The multi-input multi-output (MIMO) controller can directly adjust the Dirichlet boundary conditions on the cylinders by altering their tangential velocity.  Thus, the input operator ${B}_1$ is formed using a unit tangential velocity oriented counter-clockwise on each cylinder (forming 3 columns).

\subsubsection{Mesh Convergence Study}
    A mesh convergence study was performed with 3 additional mesh refinement levels, including one that is coarser than the mesh in figure~\ref{fig:MeshFigure} (denoted Medium) and two that are finer. The results of the study are presented in table~\ref{tab:meshConv}. After the flow is fully developed, the global energy $E(t)$ is calculated as the $L_2$ norm of the velocity field at time $t$,
\[ E(t)  \equiv \frac{1}{2}\sqrt{(\int_{\Omega} |\boldsymbol{v}(t)|^2 \operatorname{d\Omega} )}. \] This is averaged over multiple shedding periods with an averaging window time $T$. The quantities reported in table~\ref{tab:meshConv} include the time-averaged and RMS values of $E(t)$:
    \begin{displaymath}
    \overline{E} = \frac{1}{T}\int_{t_0}^{t_0+T} E(t)\, dt,\qquad
    E_{\mathrm{rms}}=\left(\frac{1}{T}\int_{t_0}^{t_0+T}(E(t)-\overline{E})^2 dt\right)^{1/2}.
    \end{displaymath}
    The dominant frequency $f_{\mathrm{peak}}$ is found from the power spectral density of the energy signal $E(t)$. After removing the mean, the Fourier transform of the signal is computed, and the dominant frequency is defined as the frequency at which the power spectral density attains its maximum.
    Relative errors 
    $\varepsilon(\overline{E})$, $\varepsilon(E_{\mathrm{rms}})$, $\varepsilon(f_{\mathrm{peak}})$ are computed with respect to the Extra-Fine mesh solution (thus those entries in the last row are zero). The relative error of a quantity $q$ computed on mesh $h$ is defined as
    \[
    \varepsilon(q_h) = \left| \frac{q_h}{q_{\mathrm{ref}}} - 1 \right|,
    \]
    where $q_{\mathrm{ref}}$ denotes the value obtained on the Extra-Fine mesh. The study showed that the mesh used for this paper is satisfactory and that the solution is a good compromise between performance and accuracy. For comparison, the papers in \cite{deng2018route,maceda2021stabilization} use a mesh with 8633 nodes and 4225 elements.  Moreover, the vorticity plot at $Re_D=50$ shown in figure~\ref{fig:VorticityPlot} is closely aligned with figure 1 from \cite{deng2018route}. 

    \begin{table}
\centering
\begin{tabular}{ccccccccc}
Mesh & Nodes & Elements & $\overline{E}$ & $E_{\mathrm{rms}}$ & $f_{\mathrm{peak}}$ &
$\varepsilon(\overline{E})$ & $\varepsilon(E_{\mathrm{rms}})$ & $\varepsilon(f_{\mathrm{peak}})$ \\
Coarse     &   4059 &   2181 & 9.436 & 0.004 & 0.100 & 0.0017 & 0.8367 & 0.4012 \\
Medium     &  14199 &   7008 & 9.448 & 0.010 & 0.167 & 0.0005 & 0.6122 & 0 \\
Fine       &  53111 &  27077 & 9.456 & 0.022 & 0.167 & 0.0003 & 0.1224 & 0 \\
Extra-Fine & 206351 & 104191 & 9.453 & 0.025 & 0.167 & 0 & 0 & 0 \\
\end{tabular}
\caption{Mesh convergence study}
\label{tab:meshConv}
\end{table}

 \subsection{Model Reduction\label{sec:MOR}}
 Model order reduction (MOR) is an important ingredient of model-based flow control and finds many applications where simulation and control of large systems are prohibitively expensive.  By reducing the complexity of these systems, MOR enables more efficient simulations and control strategies without significantly compromising accuracy. 
 There are two stages to MOR; the first is choosing a reduced basis and the second is producing a dynamical system for their coefficients.  The latter typically involves Galerkin or Petrov-Galerkin projection.  Some of the commonly used basis selection strategies include Proper Orthogonal Decomposition (POD), Balanced Truncation (BT), and Interpolatory Model Order Reduction (IMOR) based on the interpolation of transfer functions. Each of these strategies has its own advantages and is chosen to satisfy the specific requirements of the problem at hand. In the following sections, we briefly describe these strategies.  Once a reduced set of (divergence-free) basis functions for velocity is found $\{ \boldsymbol{\phi}_i \}_{i=1}^r$, we write 
 \begin{displaymath}
    \boldsymbol{v}(\cdot,t) \approx \boldsymbol{v}^r(\cdot,t) \equiv \sum_{i=1}^r \boldsymbol{\phi}_i(\cdot) a_i(t).
\end{displaymath}
Note that these basis functions are usually performed using the same approximation variables used for the flow.  Here,  these are the finite element coefficients and $\boldsymbol{x}_1 \approx \boldsymbol{\Phi}\boldsymbol{x}_r$ \eqref{eq:reducedApprox}. The columns of $\boldsymbol{\Phi}$ are the finite element coefficients of the basis functions $\boldsymbol{\phi}_i$.  Thus, an orthonormal set of basis functions has coefficients that satisfy $\boldsymbol{\Phi}^T{{M}}\boldsymbol{\Phi} = {{I}}_r$ where ${{M}}$ is the finite element mass matrix.  With this identification, we can use $\boldsymbol{x}_r$ and $\boldsymbol{a}$ interchangeably.
An ordinary differential equation for the coefficients ${\boldsymbol{x}}_1(t)$ in the form of \eqref{eq:dynamicsEq} with ${\boldsymbol{f}}({\boldsymbol{x}}_1) = {{N}}({\boldsymbol{x}}_1\otimes{\boldsymbol{x}}_1)$ can be found by projection of the Navier-Stokes equations onto the subspace spanned by the reduced basis.  The result is a reduced-order model (ROM) of the velocity where \eqref{eq:reducedApprox} and \eqref{eq:RomState} provide an approximation of its finite element coefficients.
 
 \subsubsection{POD}
 The Proper Orthogonal Decomposition was introduced to study turbulent flows by Lumley in 1967~\cite{lumley1967structure} and is one of the most well-known methods for generating reduced bases for model reduction of fluids \cite{berkooz1993proper}.  It has also been used to compute a ROM for the {\it uncontrolled} pinball problem in~\cite{deng2020low}.  The method finds a basis from a set of solution snapshots with the property that this basis is optimal in its ability to reconstruct the snapshot data.
 
Generally, the snapshot matrix ${S}$ is formed from a collection of snapshots taken at times $\{ t_1, t_2, \ldots, t_p \}$:
\begin{displaymath}
  {S} = \left[ {\boldsymbol{v}}(\cdot,t_{1}) \quad \cdots \quad {\boldsymbol{v}}(\cdot,t_{p}), \right]
\end{displaymath}
which is often referred to as the {\em input collection}.  Practically, we do not work with a matrix of vector fields, but rather the vector field in each column is replaced by coefficients from a high-dimensional ($N$) basis coming from a CFD approximation to the solution.  For example, values at cell centers, Fourier coefficients, or finite element nodal values as used in this study.
A singular value decomposition (SVD) of ${S}$ is performed as ${S}=\boldsymbol{\Phi} \mathbf{\Sigma} \mathbf{P}^T$.  The first $r$ columns of the matrix $\boldsymbol{\Phi}$ (with $r\ll N$) are used as (or typically coefficients of) a low-dimensional basis for the velocity field ${\boldsymbol{v}}$.  Note that a version of the SVD is used to ensure that the basis functions $\{ \boldsymbol{\phi}_i(\cdot) \}_{i=1}^r$ are orthonormal, for example, combining a factorization ${S}=\mathbf{Q}\mathbf{R}$ with an SVD of the resulting $\mathbf{R}$. Although a useful tool for analyzing the dynamics of the flow, and many successful applications of POD have been reported in the control setting, including~\cite{ito1998reduced,rowley2006dynamics,barbagallo2009closed-loop}, POD has a number of limitations in control problems.

In \cite{stoyanov2009reduced}, it was shown that POD alone does not necessarily accurately capture control outputs (since they are not included in their construction).   Some remedies exist; see~\cite{borggaard2016goal-oriented}, but can require an expensive optimization step.  Another issue is
that a POD model of a stable flow has been shown to lead to an unstable ROM~\cite{smith2005lowdimensional}.  This can incorrectly direct the control effort to stabilize fictitious modes. 
One solution to the problem of preservation of model stability is to use balanced truncation.

\subsubsection{Balanced Truncation}
 Balanced truncation is a technique that reduces the order of a system while preserving the essential dynamics by truncating states that have minimal contribution to the input-output behavior. It involves computing a balanced realization of the system, where states are ordered by their controllability and observability, allowing the least influential states to be truncated~\cite{mullis1976synthesis,moore81principal}.  This method has been extended to nonlinear systems in~\cite{scherpen1993balancing}, but is still limited to small problem sizes.  Recent work has tried to remedy this, e.g.,~\cite{kramer2024scalable}, but it is in the early stages.  We note that some of the computational limitations of balanced truncation have been addressed with balanced POD~\cite{rowley2005model}, which requires forward and adjoint simulations for each input and output, respectively, and finds an approximate linear balancing transformation.  This avoids assembly of the full controllability and observability Gramians.  Balanced POD provides a strategy for selecting a good set of inputs for forward and adjoint simulations.  However, the need for several forward and adjoint simulations makes it challenging to implement.
 
 \subsubsection{Interpolatory Model Order Reduction}

IMOR offers several advantages for flow control problems. Unlike snapshot-based methods such as POD, IMOR does not require expensive full-order simulations of controlled systems, and hence also avoids the challenge of selecting representative training datasets. Furthermore, IMOR inherently accounts for the control inputs and observed outputs by constructing a reduced basis that preserves the nature of the transfer function. In this work, we operate directly on the large-scale FEM operators to generate a ROM that satisfies bi-tangential Hermite interpolation conditions, ensuring that the reduced system accurately replicates the input-output frequency response of the high-fidelity Navier-Stokes dynamics.
For a linear input-output system, the goal of IMOR is to find a ROM for the input-output model of the form
\begin{align}
\label{eq:IMOR_State}
\widetilde{{E}} \dot{\tilde{\boldsymbol{x}}}(t) & =\widetilde{{A}} \widetilde{\boldsymbol{x}}(t)+\widetilde{{B}} \boldsymbol{u}(t) \\
\label{eq:IMOR_Output}
\widetilde{\boldsymbol{y}}(t) & =\widetilde{{C}} \widetilde{\boldsymbol{x}}(t)+\widetilde{{D}} \boldsymbol{u}(t),
\end{align} 
with the associated Laplace transform
$$ \widetilde{{G}}(s)=\widetilde{{C}}(s \widetilde{{E}}-\widetilde{{A}})^{-1} \widetilde{{B}}+\widetilde{{D}}. $$
In the work of \cite{gugercin2013model}, see also Theorem 2 from \cite{borggaard2015model}, it was shown that, given the full-order differential algebraic equation (DAE), interpolation points $\sigma_i \in \mathbb{C}$, and tangential direction vectors $\boldsymbol{b}_i \in \mathbb{C}^m$, $i=1, \ldots, r$, a reduced-order model of the form above matches the exact transfer function at these interpolation points and directions if basis vectors $\boldsymbol{\phi}_i$ solve
$$
\left[\begin{array}{cc}
\sigma_i {E}_{11}-{A}_{11} & {A}_{21}^T \\
{A}_{21} & 0
\end{array}\right]\left[\begin{array}{c}
\boldsymbol{\phi}_i \\
\boldsymbol{z}
\end{array}\right]=\left[\begin{array}{c}
{B}_1 \boldsymbol{b}_i \\
0
\end{array}\right]
$$
for $i=1, \ldots, r$. Construct
$$
\boldsymbol{\Phi}=\left[\boldsymbol{\phi}_1, \ldots, \boldsymbol{\phi}_r\right],
$$
then the reduced model
$$
\widetilde{{G}}(s)=\underbrace{{C} \mathbf{ \Phi}}_{\widetilde{{C}}}\left(s \underbrace{\mathbf{\Phi}^T {E}_{11} \mathbf{\Phi}}_{\widetilde{{E}}}-\underbrace{\mathbf{\Phi}^T {A}_{11} \mathbf{\Phi}}_{\widetilde{{A}}}\right)^{-1} \underbrace{\mathbf{\Phi}^T {B}_1}_{\widetilde{{B}}}+\widetilde{{D}}
$$
satisfies tangential Hermite interpolation conditions with the full model: ${G}(\sigma_i)\boldsymbol{b}_i=\widetilde{{G}}(\sigma_i)\boldsymbol{b}_i$, for $i=1,\ldots,r$. Since IMOR better captures the input-output character of the control problem, this is used in our methodology.  Note that IMOR requires $\left\{ \sigma_i, \boldsymbol{b}_i \right\}_{i=1}^r$.  There are optimal choices for these parameters; these are based on having access to $\widetilde{{G}}(s)$.  Thus, a fixed point algorithm is used that typically converges in tens of iterations~\cite{gugercin2013model}.

    \subsection{Control Problem \label{sec:polynFeedback}}
     \subsubsection{Polynomial Feedback}

This section summarizes the framework that is found in \cite{Borggaard2020}. These insights are crucial to understanding the context and background of the control methodology used in this paper. Consider a nonlinear dynamical system
\begin{equation} \label{eq:dynamicsEq}
    \dot {\boldsymbol{x}}(t) = {A} \boldsymbol{x}(t)+{B} \boldsymbol{u}(t)+\boldsymbol{f}(\boldsymbol{x}(t)),
\end{equation}
with initial condition $\boldsymbol{x}(0) = \boldsymbol{x}_0 \in \mathbb{R}^n$, and control inputs $\boldsymbol{u}(t)\in\mathbb{R}^m$, where ${A} \in \mathbb{R}^{n \times n}$ and ${B} \in \mathbb{R}^{n \times m}$ are constant matrices and $\boldsymbol{f}: \mathbb{R}^n \longrightarrow \mathbb{R}^n$ is Lipschitz continuous and satisfies $\boldsymbol{f}(0)=0$ and $\nabla_{\boldsymbol{x}} \boldsymbol{f}(0)=0$. Then the optimal control problem is to find a control $\boldsymbol{u}^*(\cdot) \in L_2(0,\infty;\mathbb{R}^m)$ that solves
$$\underset{\boldsymbol{u}}{\operatorname{min}}\ J(\boldsymbol{x},\boldsymbol{u}) = \int_0^\infty \ell\left(\boldsymbol{x}(t),\boldsymbol{u}(t)\right) dt, $$ 
with running cost $\ell({\boldsymbol{x}}(t),{\boldsymbol{u}}(t))\geq 0$ subject to equation \eqref{eq:dynamicsEq}.  To characterize the optimal control, we follow the discussion in \cite{navasca2000solution} and define the {\it value function} $\upsilon$  by \[\upsilon\left(\boldsymbol{x}_0\right) = J\left( \boldsymbol{x^*}(\cdot;\boldsymbol{x}_0), \boldsymbol{u^*}(\cdot)\right), \] where $\boldsymbol{x^*}$ and $\boldsymbol{u^*}$ are the optimal state and control input found from the initial state $\boldsymbol{x}_0$. Under the requirement that $\boldsymbol{f}$, $\ell$, and $\upsilon$ are smooth enough, the Hamilton-Jacobi-Bellman equations that describe the optimal control in feedback form $\boldsymbol{u}=\mathcal{K}(\boldsymbol{x})$ are given by
 \[\begin{aligned} 
 & 0=\frac{\partial \mathbf  \upsilon}{\partial \boldsymbol{x}}(\boldsymbol{x})({A} \boldsymbol{x}+{B} \mathcal{K}(\boldsymbol{x})+\boldsymbol{f}(\boldsymbol{x}))+\ell(\boldsymbol{x}, \mathcal{K}(\boldsymbol{x})), \\ 
 & 0=\frac{\partial \mathbf  \upsilon}{\partial \boldsymbol{x}}(\boldsymbol{x}) {B}+\frac{\partial \ell}{\partial \boldsymbol{u}}(\boldsymbol{x}, \mathcal{K}(\boldsymbol{x})). 
 \end{aligned}\]

These equations are notoriously difficult to solve since they are partial differential equations in $n$ dimensions.  The quadratic-quadratic regulator problem seeks a polynomial approximation to $\boldsymbol{u}= \mathcal{K}({\boldsymbol{x}})$ in the special case where $\boldsymbol{f}(\boldsymbol{x}) = {N} ( \boldsymbol{x} \otimes \boldsymbol{x})$, where $\mathbf N \in \mathbb{R}^{n \times n^2}$ and $\otimes$ is the usual Kronecker product.  The problem is to find a control $\boldsymbol{u}(\cdot) \in L_2(0,\infty;\mathbb{R}^m)$ that solves
$$\underset{\boldsymbol{u}}{\operatorname{min}}\ J(\boldsymbol{x},\boldsymbol{u}) = \int_0^\infty \big(\boldsymbol{x}(t)^T {Q} \boldsymbol{x}(t) + \boldsymbol{u}(t)^T  {R} \boldsymbol{u}(t)
\big) dt, $$ where ${{Q}}\geq 0$ and ${{R}} > 0$ are symmetric weighing matrices for the state and the control vectors, respectively. 
 We now follow the steps shown in \cite{Borggaard2020} to derive the QQR. The value function is expanded as
\begin{equation}
\label{eq:valueExpansion}
\upsilon(\boldsymbol{x}) = \boldsymbol{\upsilon}_2^T(\boldsymbol{x} \otimes \boldsymbol{x}) + \boldsymbol{\upsilon}_3^T (\boldsymbol{x} \otimes \boldsymbol{x} \otimes \boldsymbol{x}) + \ldots 
\end{equation}
where the coefficients are $\boldsymbol{\upsilon}_q \in \mathbb{R}^{n^q}$. The feedback operator is expanded similarly as
\begin{equation} \label{eq:feedbackExpansion}
\boldsymbol{u}(\boldsymbol{x}) = \boldsymbol{k}_1 \boldsymbol{x} + \boldsymbol{k}_2 (\boldsymbol{x} \otimes \boldsymbol{x}) + \boldsymbol{k}_3 (\boldsymbol{x} \otimes \boldsymbol{x} \otimes \boldsymbol{x}) + \ldots 
\end{equation}
where $\boldsymbol{k}_q \in \mathbb{R}^{m \times n^q}$. The QQR derivation proceeds to find the value of feedback operator approximations to any given degree by matching coefficients of lowest degree first and sequentially finding the next degree terms by solving Kronecker sum systems. The lowest degree coefficients $\boldsymbol{\upsilon}_2$ solve the algebraic Riccati equation and $\boldsymbol{k}_1$ is the gain associated with the Linear Quadratic Regulator (LQR). Higher degree coefficients are found using linear systems with special Kronecker sum systems defined for an $n\times n$ matrix $\mathbf{X}$ as
$$
\mathcal{L}_d(\mathbf{X}) \equiv \underbrace{\mathbf{X} \otimes \cdots \otimes \mathbf{I}_n}_{d \text { terms }}+\cdots+\underbrace{\mathbf{I}_n \otimes \cdots \otimes \mathbf{X}}_{d \text { terms }}
$$
where $\mathbf{I}_n$ is the $n\times n$ identity matrix. The next term in the value function approximation can be found by solving
\begin{equation} \label{eq:higherValueFuncEq}
\mathcal{L}_3\left({A} + {B} \boldsymbol{k}_1^T\right) \boldsymbol{\upsilon}_3=-\mathcal{L}_2\left({N}^T\right) \boldsymbol{\upsilon}_2. 
\end{equation}
Higher degree value function approximations can be found using similar formulas. For more details, see \cite{Borggaard2020}. The formula for the next term in the approximation of the controller can be obtained from
\begin{equation} \label{eq:gainEq}
\boldsymbol{k}_d=- {R}^{-1}\left(\mathcal{L}_{d+1}\left({B}^{T}\right) \boldsymbol{\upsilon}_{d+1}\right)^{T}.
\end{equation}
Note that systems of the form \eqref{eq:higherValueFuncEq} have $n^d$ unknowns and become prohibitively expensive for large values of $d$ and $n$.  Therefore, we truncate our approximations of \eqref{eq:valueExpansion} and \eqref{eq:feedbackExpansion} to the terms $\boldsymbol{\upsilon}_3\in \mathbb{R}^{n^3}$ and ${\boldsymbol{k}}_2\in \mathbb{R}^{m\times n^2}$ in this study.

\subsubsection{QQR for the Pinball Problem\label{sec:QQR}}
To build a quadratic model for the pinball problem, we built the reduced basis $\boldsymbol{\Phi}$ using IMOR then used it to project the finite element model in \eqref{eq:CFD_DAE} to
find
\begin{equation} \label{eq:quadraticROM}
  \widetilde{{E}} \dot{\tilde{\boldsymbol{x}}}(t)  =\widetilde{{A}} \widetilde{\boldsymbol{x}}(t)+\widetilde{{B}} \boldsymbol{u}(t) + \widetilde{{N}}\left( \widetilde{\boldsymbol{x}} \otimes \widetilde{\boldsymbol{x}} \right).
\end{equation}
Since $\widetilde{{E}}$ is invertible, we exactly match the quadratic-quadratic regulator problem defined in section \ref{sec:polynFeedback}, given the dynamics in \eqref{eq:quadraticROM}, we seek a control $\boldsymbol{u}(\cdot) \in L_2(0,\infty;\mathbb{R}^m)$ that solves
\begin{align*}
\underset{\boldsymbol{u}}{\operatorname{min}}\ J_r(\tilde{\boldsymbol{x}},\boldsymbol{u}) &= \int_0^\infty \big(\tilde{\boldsymbol{y}}(t)^T  \tilde{\boldsymbol{y}}(t) + \boldsymbol{u}(t)^T  {R} \boldsymbol{u}(t)
\big) dt \\
&= \int_0^\infty \big(  \tilde{\boldsymbol{x}}(t)^T \tilde{{C}}^T \tilde{{C}} \tilde{\boldsymbol{x}}(t) + \boldsymbol{u}(t)^T  {R} \boldsymbol{u}(t)
\big) dt,
\end{align*} where $ {R} = {I}_3$ for the pinball problem. 
Since
\[ \tilde{\boldsymbol{x}}^T(t)\tilde{{C}}^T \tilde{{C}} \tilde{\boldsymbol{x}}(t) = \|\tilde{{C}} \tilde{\boldsymbol{x}}(t)\|^2 = \| {C}\boldsymbol{\Phi}\tilde{\boldsymbol{x}}(t)\|^2 \approx \| {C}\boldsymbol{x}_1(t)\|^2, \]
a successful control should drive averages of $\boldsymbol{x}_1$ in the patches of figure~\ref{fig:Cmatrix} to 0 over time.  As $\boldsymbol{x}_1$ are finite element coefficients of $\boldsymbol{v} - \boldsymbol{v_{ss}}$ it is reasonable to expect the field $\boldsymbol{v} \rightarrow \boldsymbol{v}_{ss}$.  It is important to make principled choices about where the controlled outputs are located.  Recall that we developed ${C}$ based on an approximation of the vortex formation length.

\subsection{Closing the Loop\label{sec:closedLoop}}
The ROM is utilized to design both the linear and the QQR controllers. For the QQR controller, the nonlinear term is incorporated in the control design, as seen in \eqref{eq:quadraticROM}.

The linear controller matches the LQR solution when ignoring the quadratic term and is consistent with $\boldsymbol{k}_1$ from \eqref{eq:feedbackExpansion}.  We will denote this as $\boldsymbol{u}^{(1)}(t) = \boldsymbol{k}_1\tilde{\boldsymbol{x}}$, where $\boldsymbol{k}_1$ can be derived from the solution of the algebraic Riccati equation
\[\widetilde{{A}}^T \mathcal{V}_2 \widetilde{{E}}+\widetilde{{E}}^T \mathcal{V}_2 \widetilde{{A}}-\widetilde{{E}}^T \mathcal{V}_2 \widetilde{{B}}  {R}^{-1} \widetilde{{B}}^T \mathcal{V}_2 \widetilde{{E}}+\widetilde{{C}}^T \widetilde{{C}}=0,\] where $\mathcal {V}_2$ is the quadratic value function approximation term, $\boldsymbol{\upsilon}_2$, reshaped into a square matrix.  This is computed as $\boldsymbol{k}_1 = -{R}^{-1}\widetilde{{B}}^T\mathcal{V}_2$.  Regarding the quadratic controller, starting with $\boldsymbol{\upsilon}_2$ from the LQR, $\boldsymbol{\upsilon}_3$ can be solved for using \eqref{eq:higherValueFuncEq}. Then, we can solve for $\boldsymbol{k}_2$ using \eqref{eq:gainEq}.
This is the QQR algorithm applied to the reduced order model \eqref{eq:quadraticROM}. Using the same notation, we will denote the quadratic controller as $\boldsymbol{u}^{(2)}(t) = \boldsymbol{k}_1\tilde{\boldsymbol{x}}(t) + \boldsymbol{k}_2\left( \tilde{\boldsymbol{x}}(t)\otimes \tilde{\boldsymbol{x}}(t)\right)$.  

It is important to emphasize that in this work, the controller is being evaluated in the full-order system \eqref{eq:CFD_DAE}.  So interpreting the controller action with the $\boldsymbol{x}_1$ variables is helpful.  We can lift the gains to the full-order model dimension combining the two calculations.  If we interpret $\tilde{\boldsymbol{x}}_i(t)$ as the finite element approximation of $\int_\Omega \boldsymbol{\phi}_i(x) \boldsymbol{v}'(x,t)\ dx$, with $\boldsymbol{\phi}_i$ and $\boldsymbol{v}'$ written in the same finite element basis, this leads to \(\tilde{\boldsymbol{x}}=\mathbf{\Phi}^T {M} \boldsymbol{x}_1\) where ${M}$ is the finite-element mass matrix and the control law of interest.  For example, the lifted full-order feedback law of degree 2 is
\begin{align*}
\boldsymbol{u}^{(2)} = \mathbf{\mathcal{K}}(\tilde{\boldsymbol{x}}) &= \boldsymbol{k}_1 \tilde{\boldsymbol{x}} + \boldsymbol{k}_2(\tilde{\boldsymbol{x}}\otimes \tilde{\boldsymbol{x}}) \\
&= \boldsymbol{k}_1 \mathbf{\Phi}^T {M} \boldsymbol{x}_1 + \boldsymbol{k}_2(\mathbf{\Phi}^T {M} \boldsymbol{x}_1\otimes \mathbf{\Phi}^T {M} \boldsymbol{x}_1)\\
&\equiv \bar{\boldsymbol{k}}_1  \boldsymbol{x}_1 + \bar{\boldsymbol{k}}_2( \boldsymbol{x}_1\otimes \boldsymbol{x}_1).
\end{align*}
Here, $\bar{\boldsymbol{k}}_1 = \boldsymbol{k}_1 \mathbf{\Phi}^T {M}$ and $\bar{\boldsymbol{k}}_2 = \boldsymbol{k}_2(\mathbf{\Phi}^T {M}\otimes \mathbf{\Phi}^T {M})$ and can be precomputed.

\section{Numerical Results}

\begin{figure}[!h]
    \centering
    \begin{subfigure}{0.48\textwidth}
        \centering
       \includegraphics[width=\textwidth]{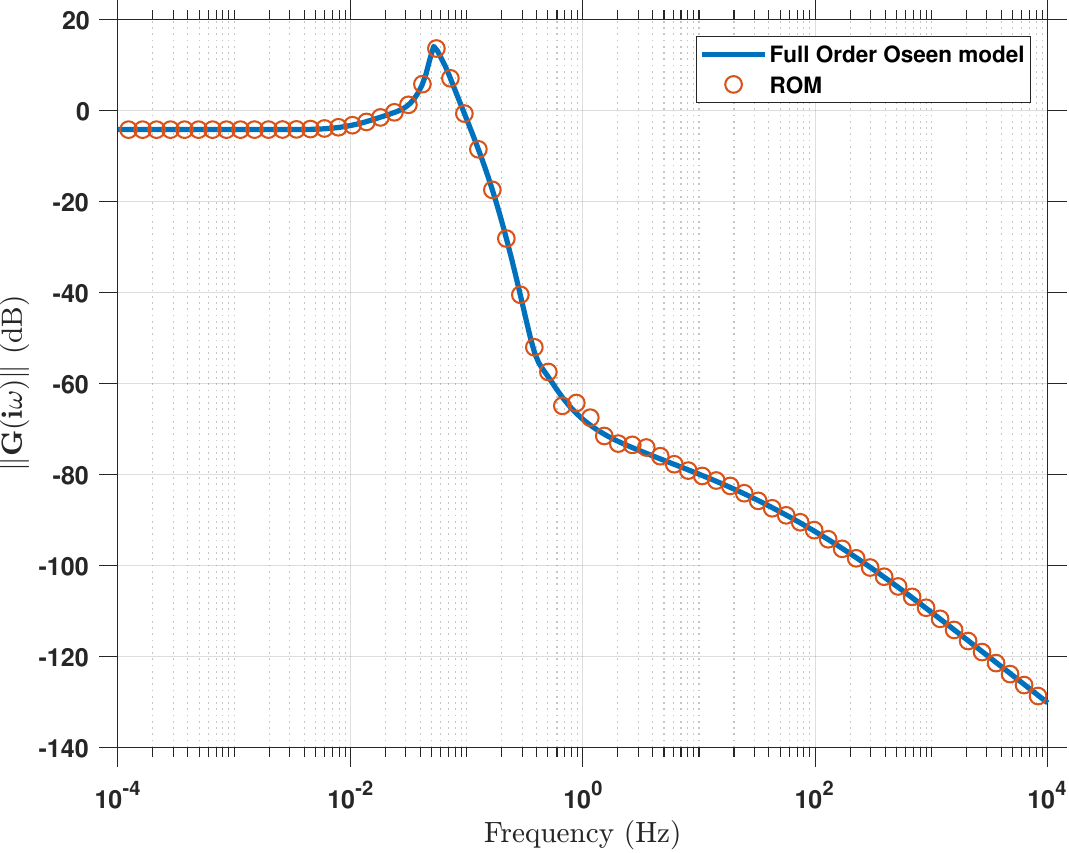} 
    \caption{} \label{fig:BodePlot}
    \end{subfigure}
    \hfill
    \begin{subfigure}{0.48\textwidth}
            \centering
             \includegraphics[width=\textwidth]{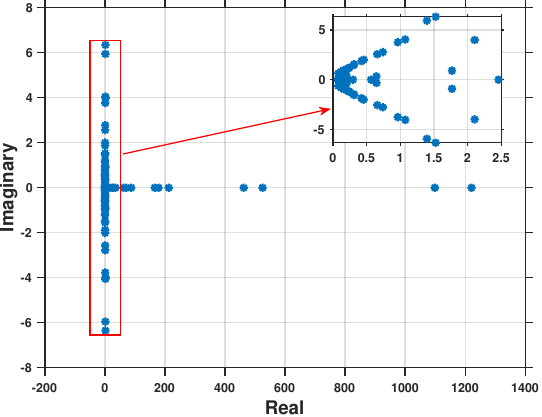} 
            \caption{} \label{fig:IRKAShifts}
    \end{subfigure}
    \caption{(a) Bode plot of the IRKA ROM model with $r=80$ compared with the full-order model for $Re_D=30$. (b) Final set of shifts $\{\sigma_i\}_{i=1}^r$, found by IRKA for $Re_D=30$ .}
\end{figure}

In this section, we present the numerical results obtained from simulations conducted with two different Reynolds numbers, \(Re_D = 30\) and \(Re_D = 50\). In figure \ref{fig:BodePlot}, the Bode plot of the ROM at $Re_D=30$ shows good agreement with the FOM, and indicates that the interpolation conditions are approximately satisfied.  Similar results are found for the \(Re_D = 50\) case, but are not shown here for brevity.  Throughout the updates of the IRKA algorithm, the tangential directions $\boldsymbol{b}_i$ are updated from the current iterate of the ROM by solving a generalized eigenvalue problem with $(\widetilde{{A}}^T,\widetilde{{E}}^T)$, then setting $\boldsymbol{b}_i = \widetilde{{B}}^T \boldsymbol{w}_i$, where $\boldsymbol{w}_i$ is the $i\,$th eigenvector. In each iteration of the IRKA algorithm, the $i$th eigenvalue, $\lambda_i$, is used to compute the $i\,$th shift as $\sigma_i = -\lambda_i$. The final shifts $\{\sigma_i\}_{i=1}^r$, found by IRKA at $Re_D=30$, are shown in figure \ref{fig:IRKAShifts}. 

In the rest of the section, the performance of the QQR controller is compared with that of a linear controller to evaluate their effectiveness in driving the flow towards a steady-state solution. For \(Re_D = 30\), we analyze the \(L_2\)-norm error in velocity, control input, and the target running cost. 
In addition, we examine the lift and drag coefficients to assess the impact of control on the body forces of each of the cylinders. Similarly, for \(Re_D = 50\), we explore the ability of the control  to achieve the desired steady-state solution and its influence on error metrics and cost functions. The subsequent subsections provide a detailed comparison of the performance of the controllers, highlighting the advantages of the QQR controller in terms of faster convergence and lower running costs.

\subsection{Case with $Re_D=30$}
\begin{figure}[!htb]
    \centering
    \begin{subfigure}{.4\columnwidth} 
        \framedgraphic{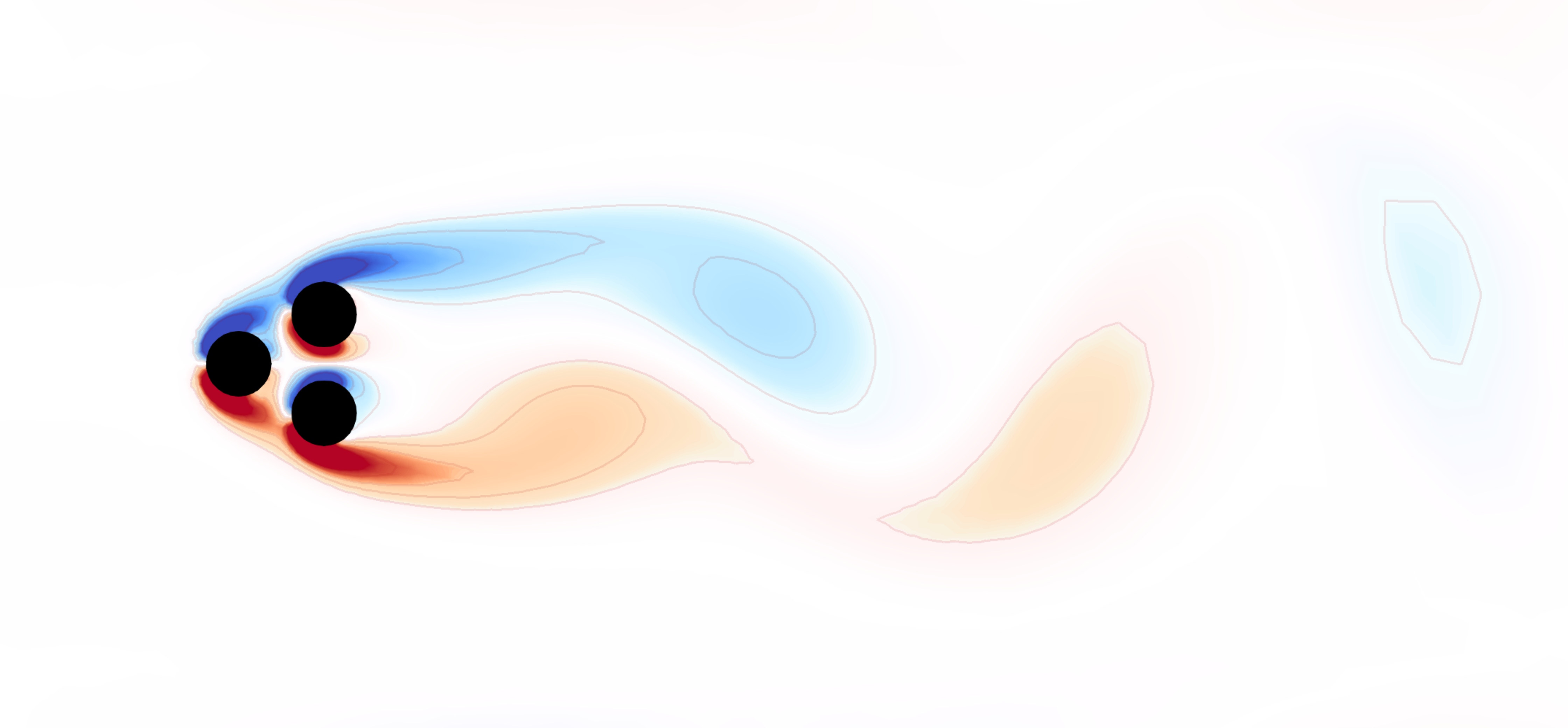}
        \label{fig:Re30Vort_t=10}
    \end{subfigure}
    \hspace{0.1cm}
    \begin{subfigure}{.4\columnwidth}
        \framedgraphic{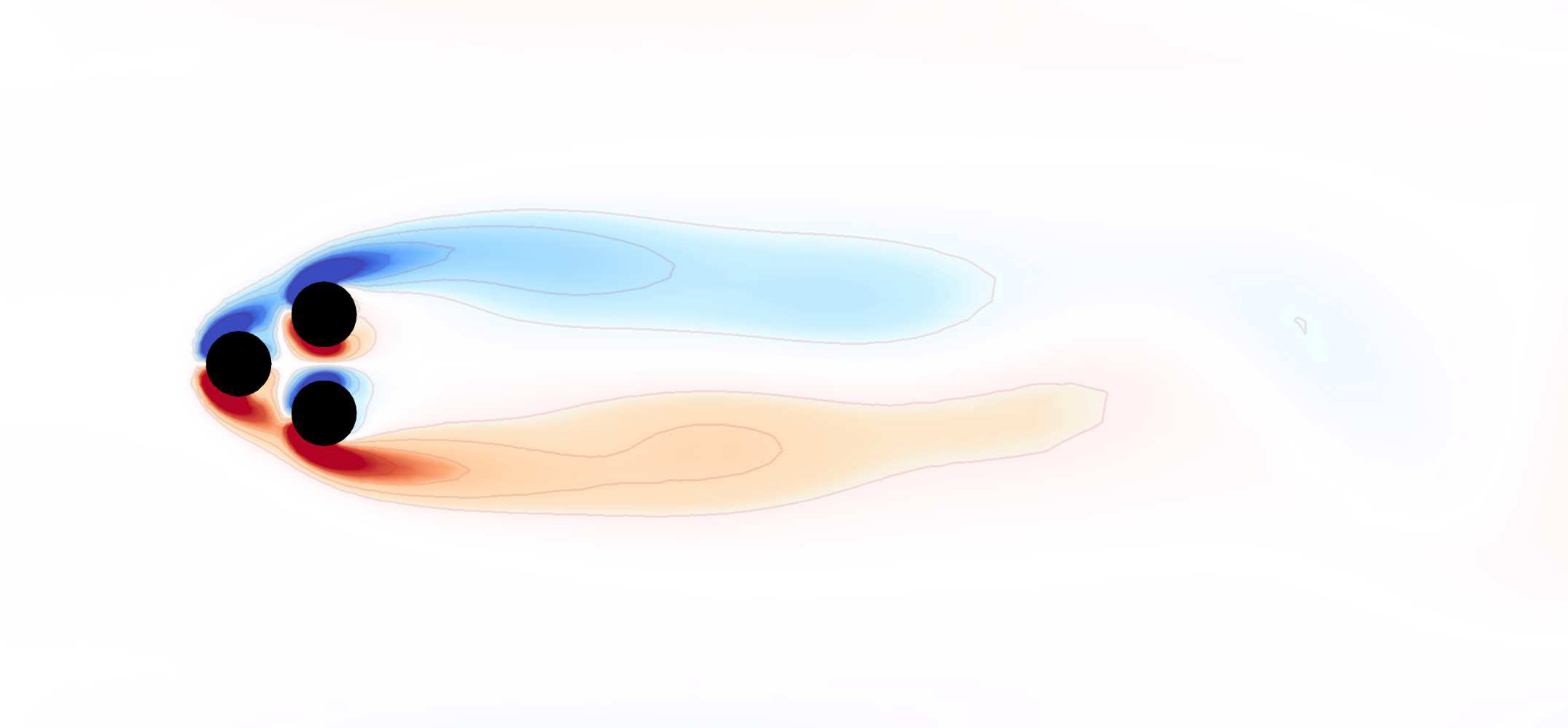}
        \label{fig:Re30Vort_t=30}
    \end{subfigure}\\[1ex]
    \begin{subfigure}{.4\columnwidth}
        \centering
        \framedgraphic{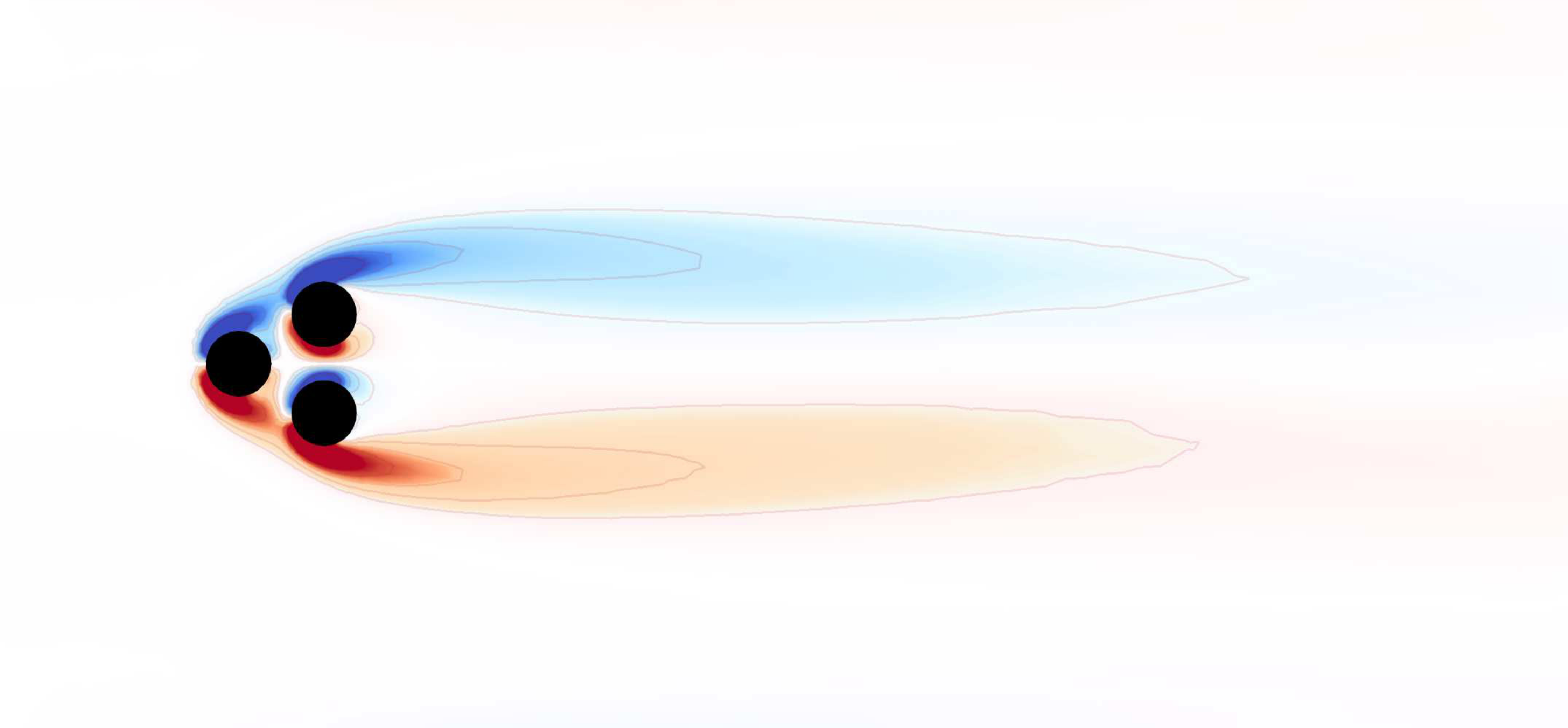}
        \label{fig:Re30Vort_t=80}
    \end{subfigure}
    \caption{Vorticity fields for $Re_D=30$ with QQR controller $\boldsymbol{u}^{(2)}$ at three different times; \subref{fig:Re30Vort_t=10} t = 10, \subref{fig:Re30Vort_t=30} t = 30, \subref{fig:Re30Vort_t=80} t = 80.}
    \label{fig:Re30Vorticity}
\end{figure}

\begin{figure}[H]
    \centering
    \begin{subfigure}{.4\columnwidth} 
    \includegraphics[width=\linewidth]{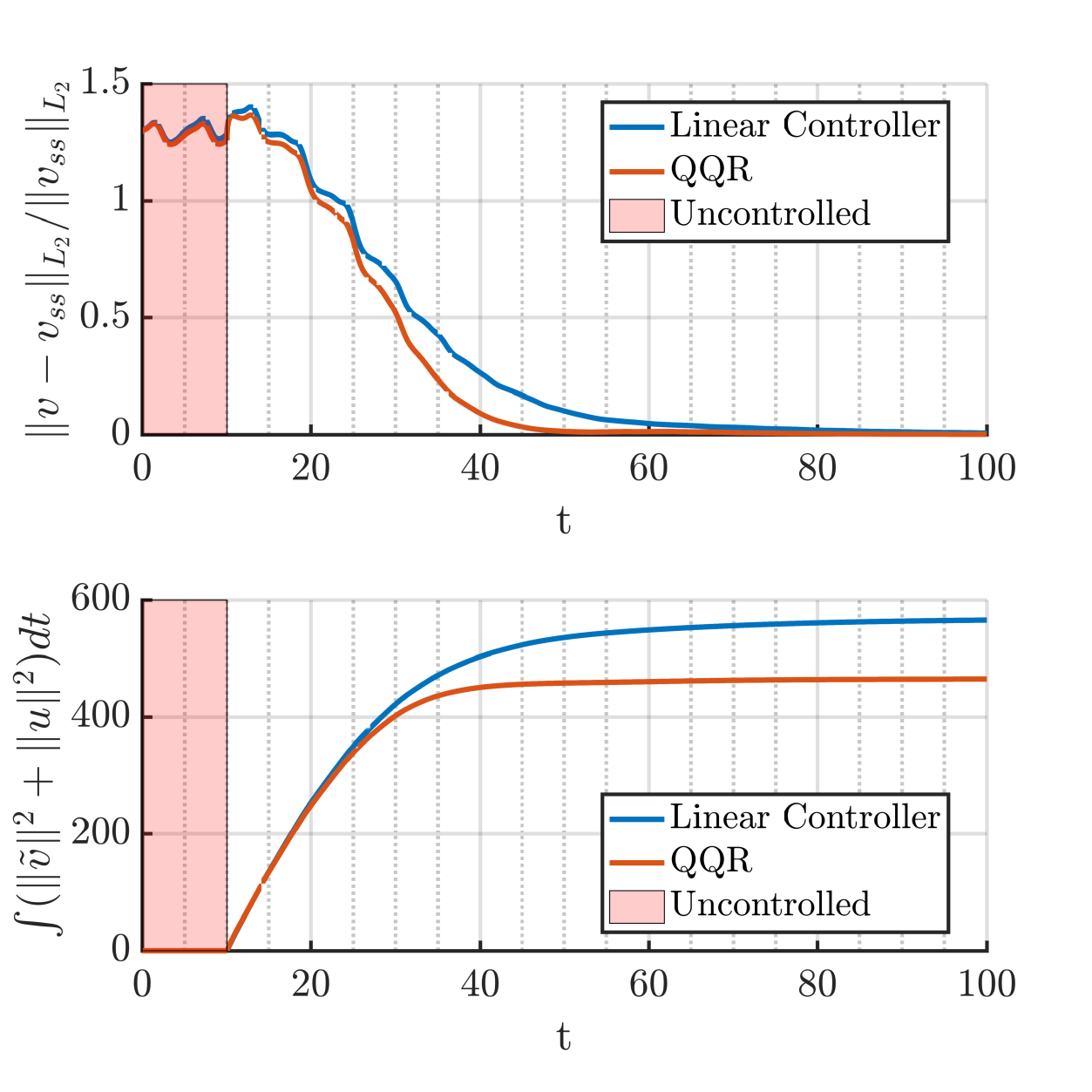} 
    \caption{} \label{fig:Re30Error}
    \end{subfigure}\vspace{-.05in}
    \begin{subfigure}{.4\columnwidth}
    \includegraphics[width=\linewidth]{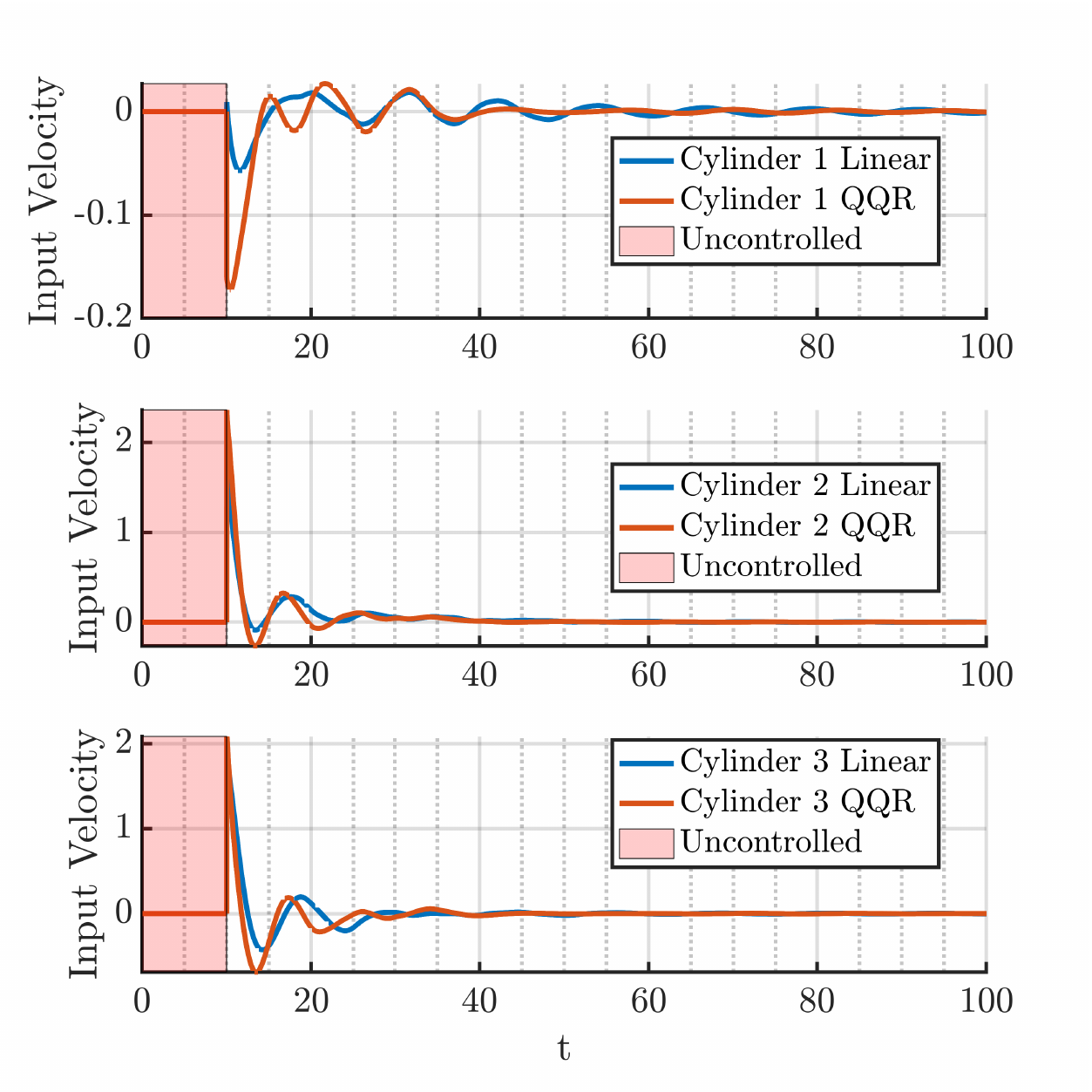}
    \caption{} \label{fig:Re30Control}
    \end{subfigure}\vspace{-.05in}
    \caption{$Re_D=30$ case \subref{fig:Re30Error} The top figure shows a comparison between the $L_2$ error in velocity for the linear and QQR controllers, while the bottom one shows a comparison of the integrated running cost for both controllers. \subref{fig:Re30Control} Shows the control efforts for the three cylinders and both controllers, represented by their tangential velocity.} 
\end{figure}

The vorticity field of the flow, shown in figure \ref{fig:Re30Vorticity}, is plotted at $t=10, 30,$ and $100$s. In figure \ref{fig:Re30Error}, we plot the ${L}_2$ norm of the velocity perturbation over the whole domain, which is the difference between the instantaneous flow field and the steady-state flow field. The commanded input, which is the tangential velocity of the three cylinders, is seen in figure \ref{fig:Re30Control}. It can be noted that the QQR control drives the flow to the steady-state within a relative error of less than 2\%,  about 40.1\% faster than the LQR control. The norm of the steady-state velocity was found to be $\|\boldsymbol{v}_{ss}\|_{L_2} = 19.1$. Moreover, the ideal running cost function plot, shows that the total cost of the QQR controller, which combines the control effort and the state error costs, is 13.5\% lower than the cost of the linear controller.

\begin{figure}[!h]
    \centering
    \begin{subfigure}{.4\columnwidth} 
    \includegraphics[width=\linewidth]{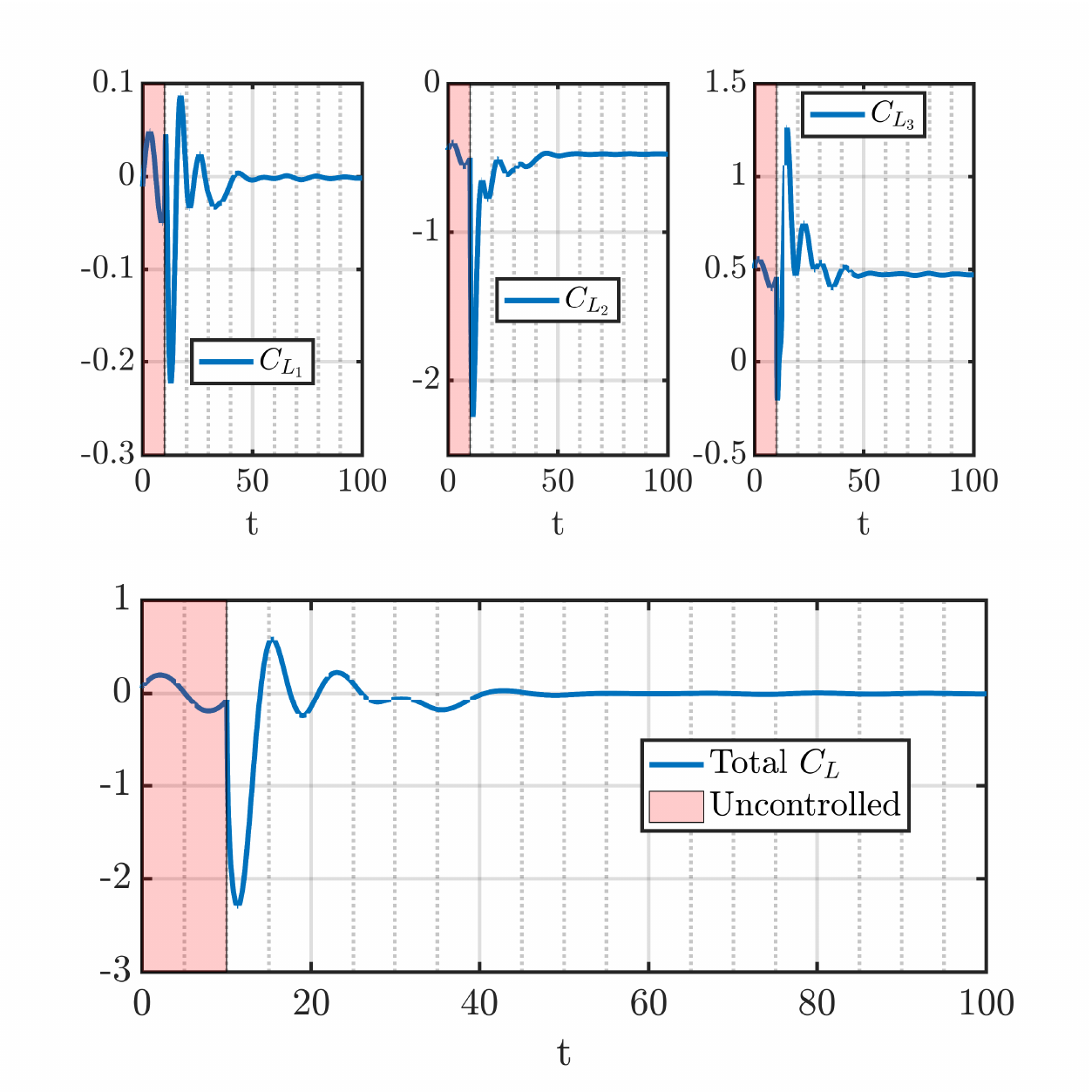} 
    \caption{} \label{fig:Re30Lift}
    \end{subfigure}\vspace{-.05in}
    \begin{subfigure}{.4\columnwidth}
    \includegraphics[width=\linewidth]{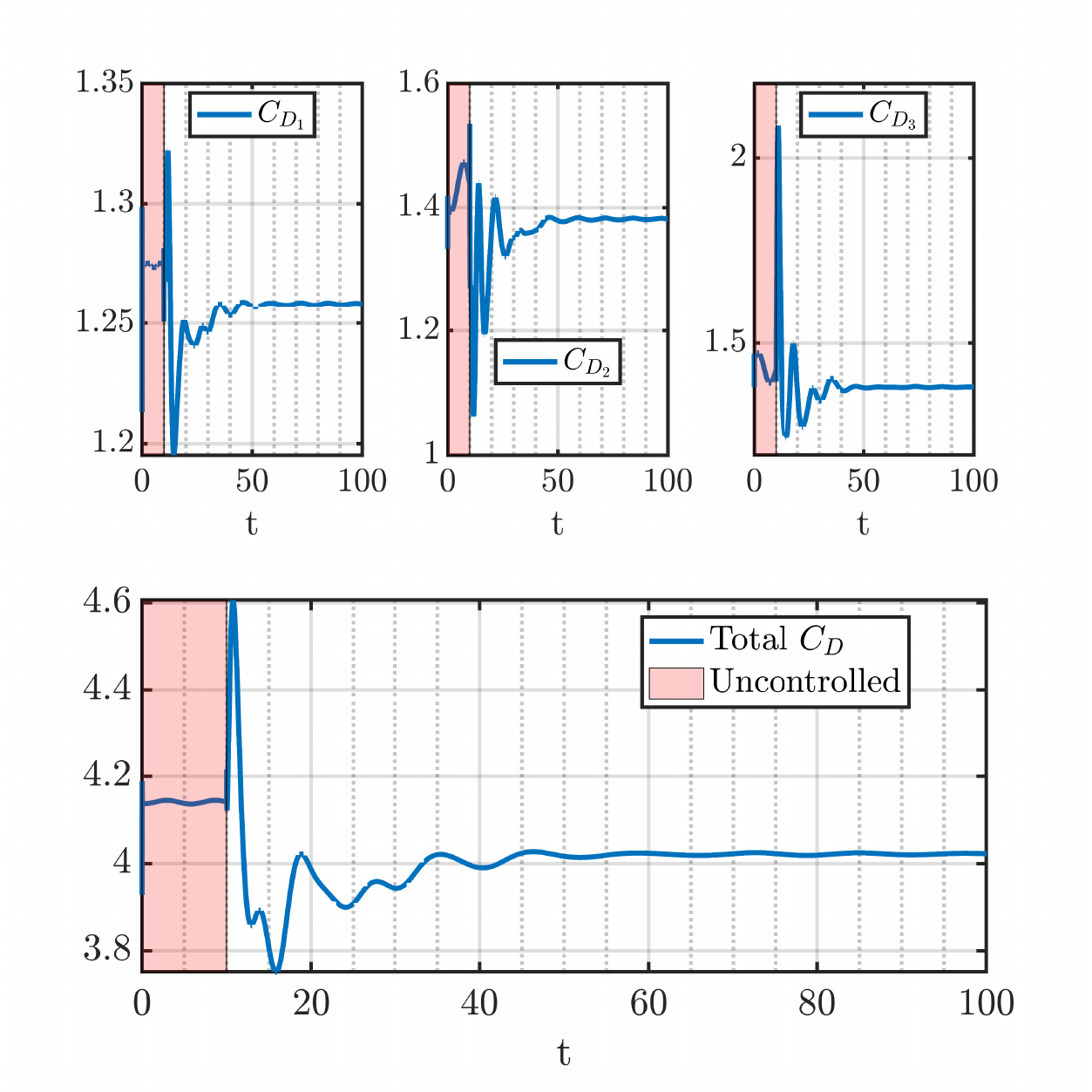}
    \caption{} \label{fig:Re30Drag}
    \end{subfigure}\vspace{-.05in}
    \caption{$Re_D=30$ case \subref{fig:Re30Error} The top figure presents the lift coefficients on each individual cylinder while the bottom figure provides the evolution of $C_L$ on the system of cylinders. \subref{fig:Re30Control} Figures corresponding to the drag coefficient $C_D$ evolution.} 
\end{figure}

In figures \ref{fig:Re30Lift} and \ref{fig:Re30Drag}, the lift and drag coefficients are plotted, respectively. Although not directly the goal of this paper, the drag was reduced by approximately 2.9\% (both for linear and QQR controllers), while vibration-inducing oscillations in lift were reduced from a magnitude of 0.198 to 0.007, after about 50 time units of actuation. Note that the lift approaches zero as time goes to infinity.

\subsection{Case with $Re_D = 50$}
\begin{figure}[H]
    \centering
    \begin{subfigure}{.4\columnwidth} 
            \framedgraphic{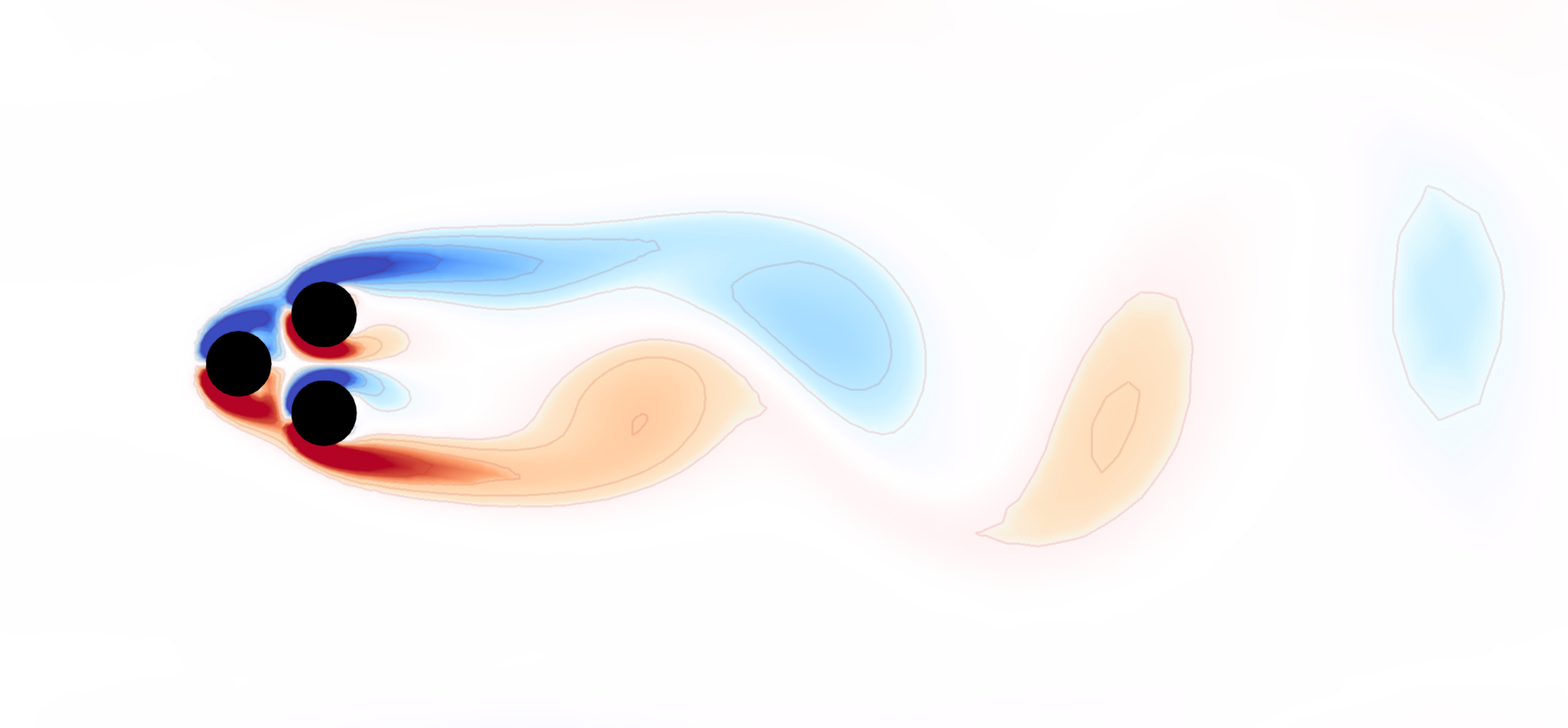}
        \label{fig:Re50Vort_t=10}
    \end{subfigure}
    \hspace{0.1cm}
    \begin{subfigure}{.4\columnwidth}
        \framedgraphic{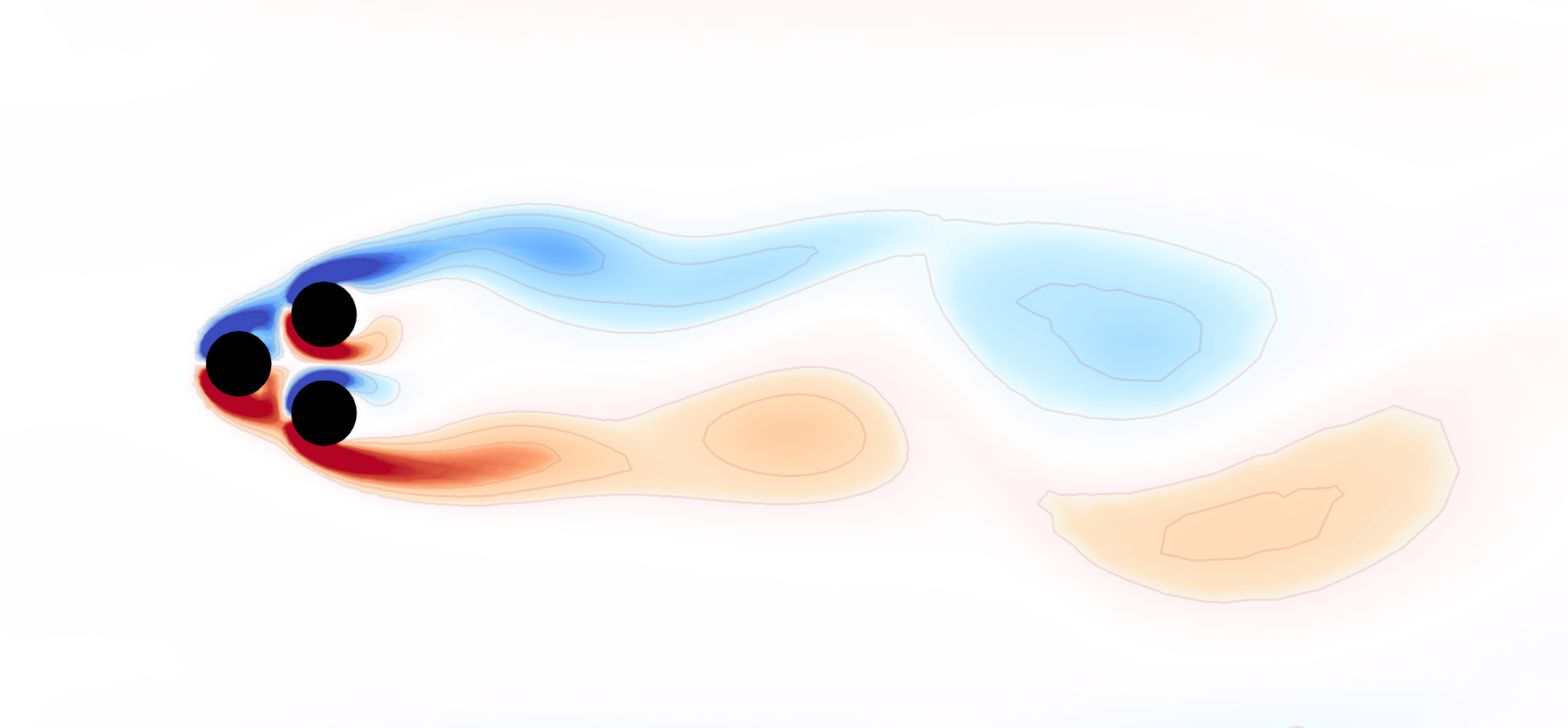}
        \label{fig:Re50Vort_t=30}
    \end{subfigure}\\[1ex]
    \begin{subfigure}{.4\columnwidth}
        \framedgraphic{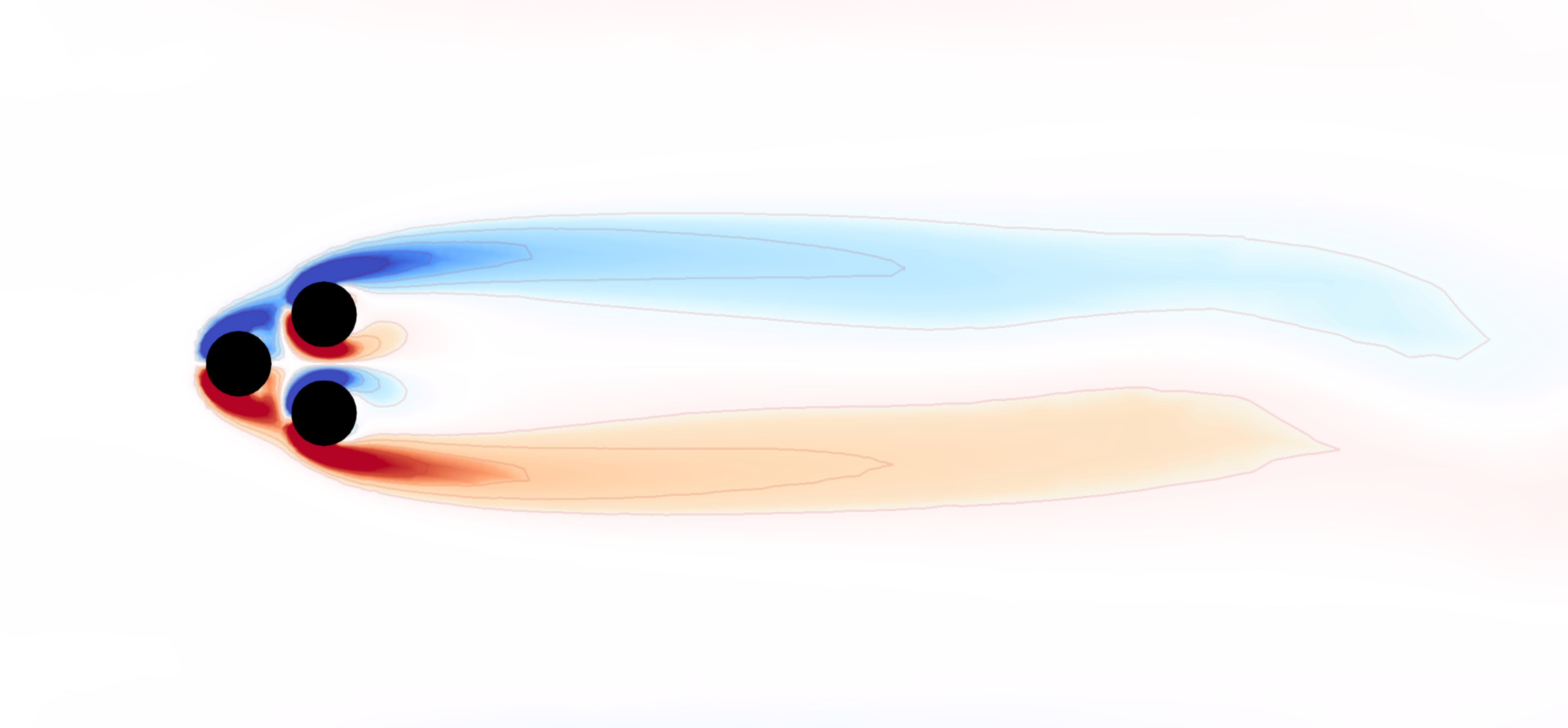}
        \label{fig:Re50Vort_t=80}
    \end{subfigure}
    \caption{Vorticity fields for $Re_D=50$ with the QQR controller $\boldsymbol{u}^{(2)}$ at three different times; \subref{fig:Re50Vort_t=10} t = 10, \subref{fig:Re50Vort_t=30} t = 30, \subref{fig:Re50Vort_t=80} t = 80.}
    \label{fig:Re50Vorticity}
\end{figure}

\begin{figure}[H]
    \centering
    \begin{subfigure}{.4\columnwidth}
    \includegraphics[width=\linewidth]{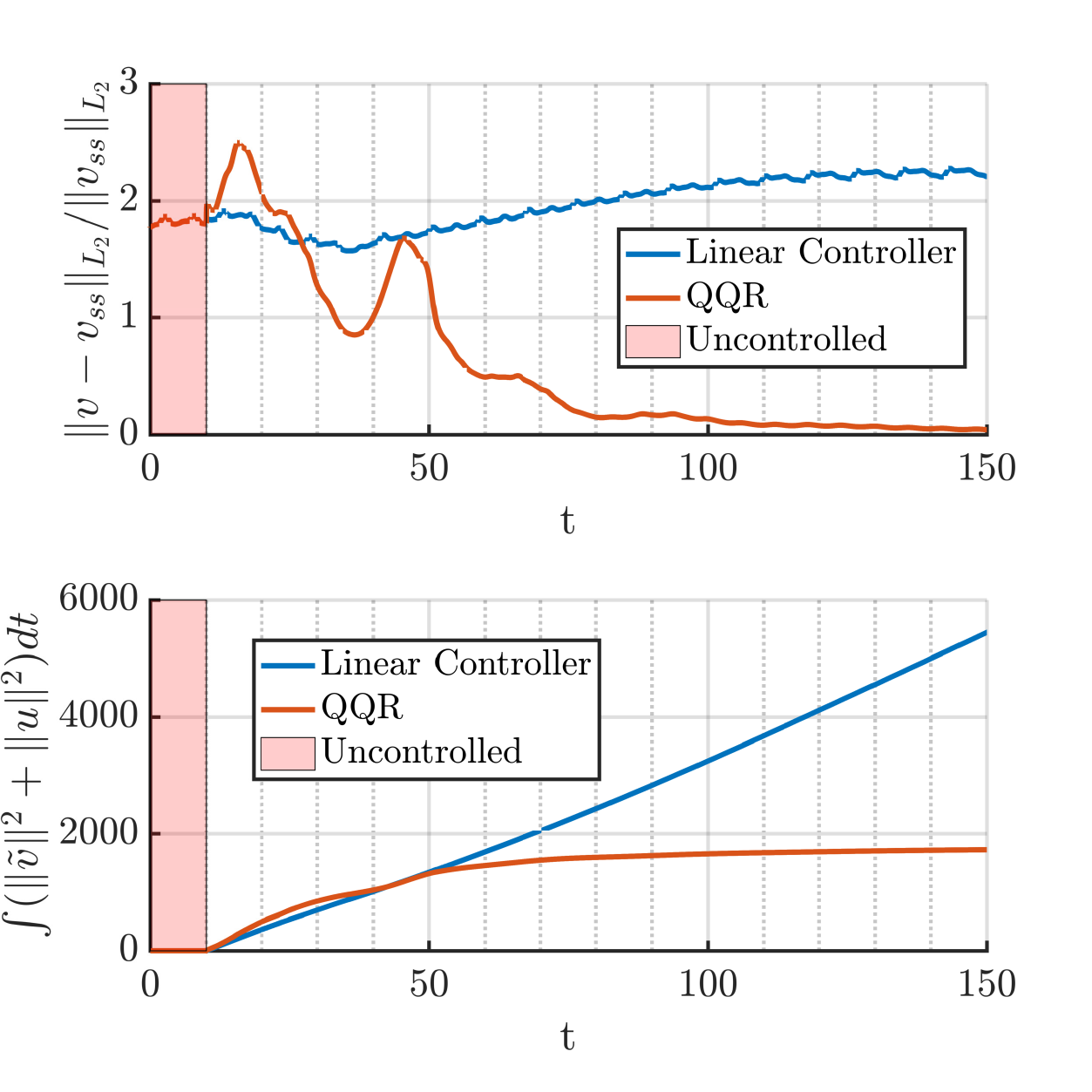} 
    \caption{}
    \label{fig:Re50Error}
    \end{subfigure}\vspace{-.05in}
    \begin{subfigure}{.4\columnwidth}
    \includegraphics[width=\linewidth]{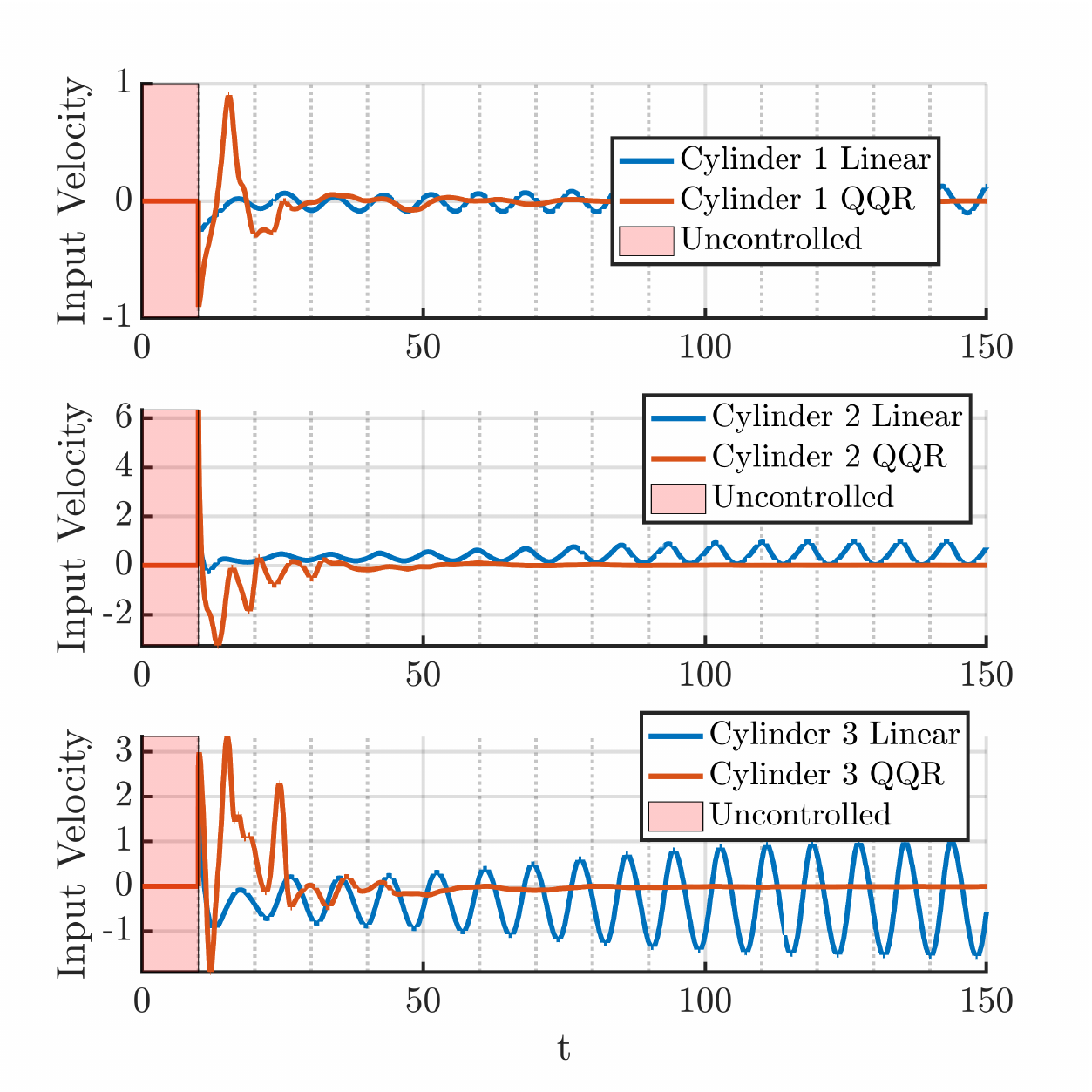}
    \caption{}
    \label{fig:Re50Control}
    \end{subfigure}\vspace{-.05in}
    \caption{$Re_D=50$ case \subref{fig:Re50Error} The top figure shows a comparison between the $L_2$ error in velocity for the linear and QQR controllers, while the bottom one shows a comparison between the ideal cost function for both controllers. \subref{fig:Re50Control} Shows the control efforts for the three cylinders and both controllers, represented by their tangential velocity}

\end{figure}

In this case, the linear controller cannot drive the flow to the desired steady-state solution. Instead, it drives it into a different limit cycle that is further away from the desired steady-state solution. However, the QQR can control the flow and drive it to the desired solution. The vorticity plot for this case is shown in figure \ref{fig:Re50Vorticity}. The $L_2$ error in velocity deviation and the ideal running cost function are shown in figure \ref{fig:Re50Error}, and the tangential commanded velocities for the three cylinders are shown in figure \ref{fig:Re50Control}. The norm of the steady-state velocity was found to be $ \|\boldsymbol{v}_{ss}\|_{L_2} =  19.3$.

\begin{figure}[H]
    \centering
    \begin{subfigure}{.45\columnwidth}
    \includegraphics[width=\linewidth]{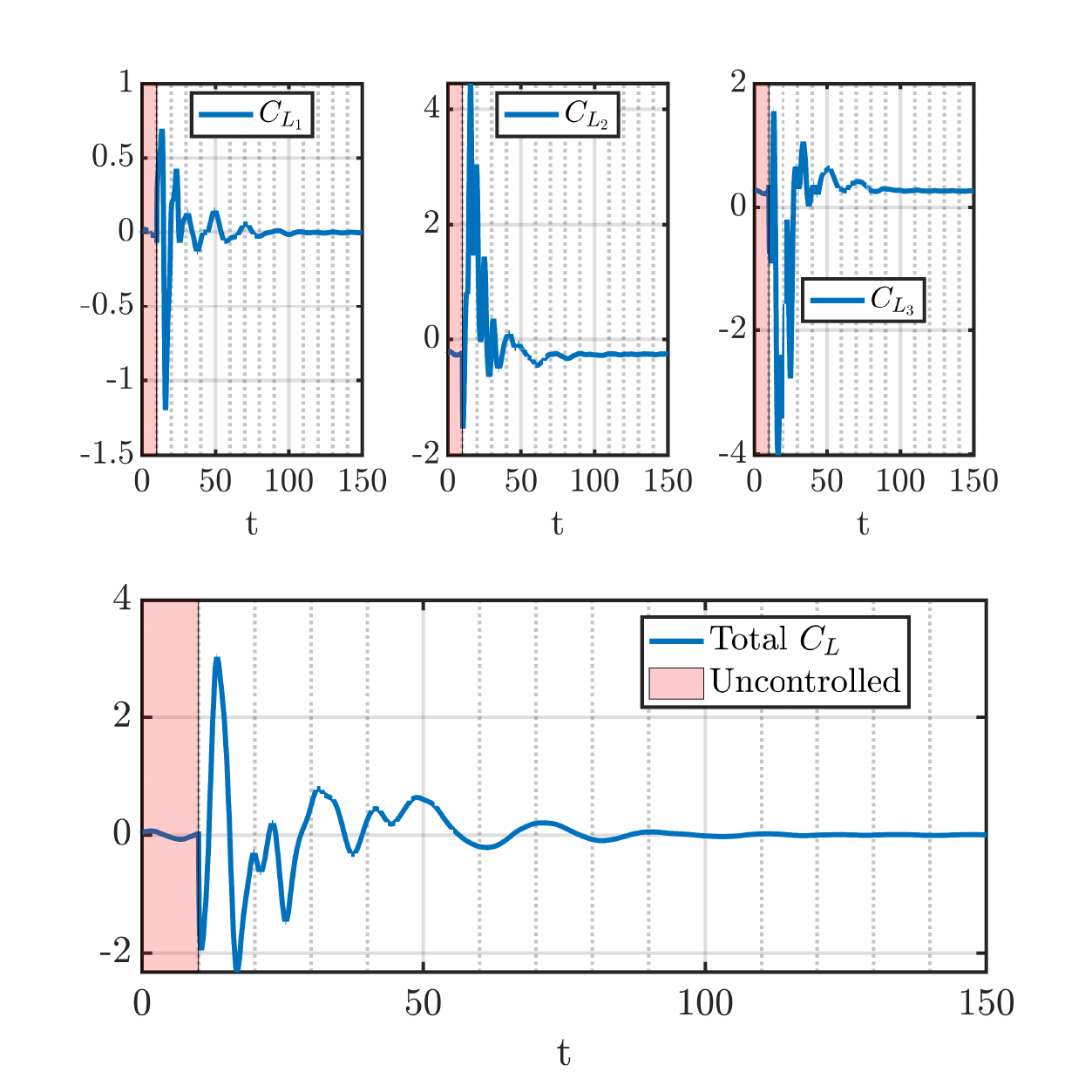} 
    \caption{}
    \label{fig:Re50Lift}
    \end{subfigure}\vspace{-.05in}
    \begin{subfigure}{.45\columnwidth}
    \includegraphics[width=\linewidth]{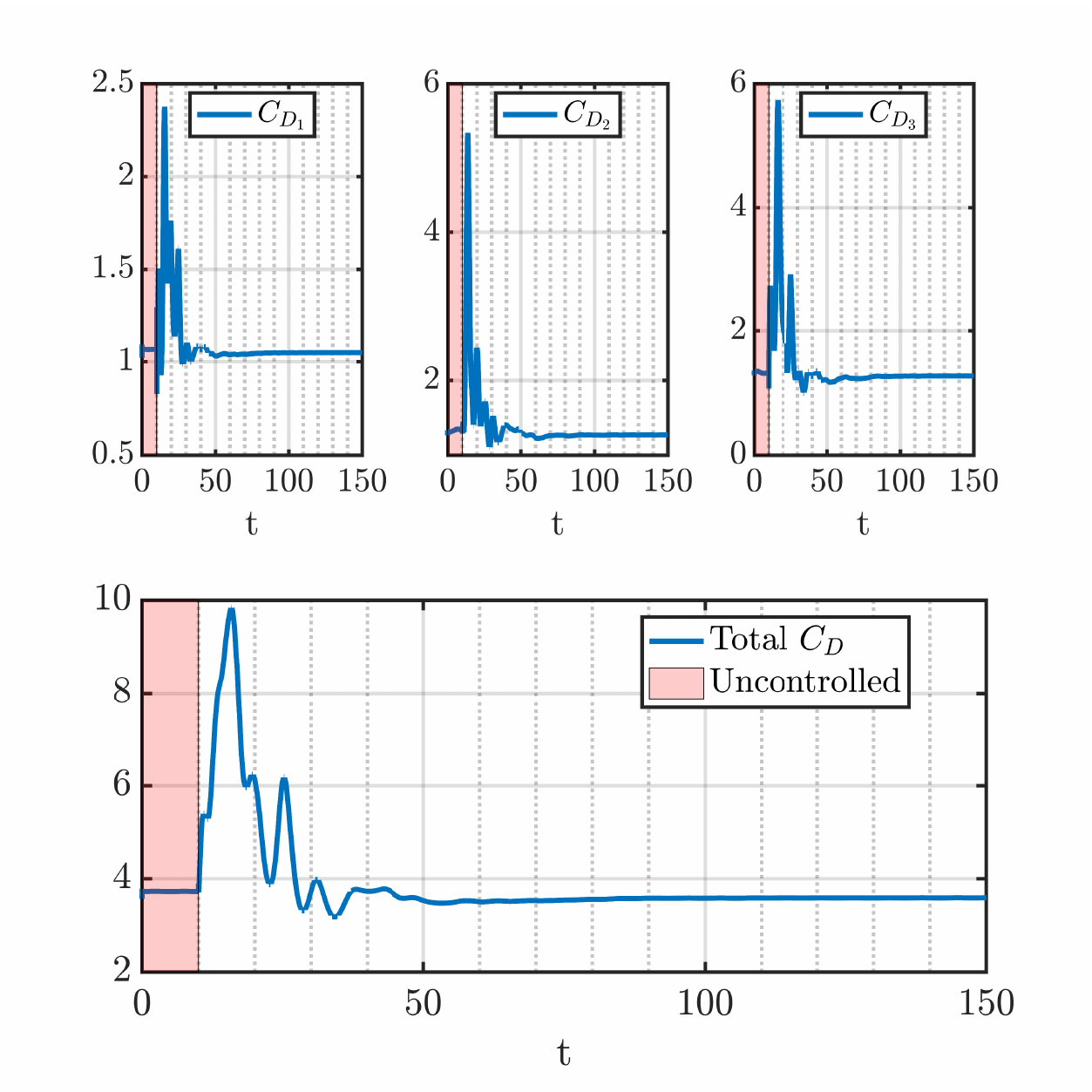}
    \caption{}
    \label{fig:Re50Drag}
    \end{subfigure}\vspace{-.05in}
    \caption{$Re_D=50$ case \subref{fig:Re30Error} The top figure presents the lift coefficients on each individual cylinder while the bottom figure provides the evolution of $C_L$ on the system of cylinders. \subref{fig:Re30Control} Figures corresponding to the drag coefficient $C_D$ evolution.} 

\end{figure}

In figures \ref{fig:Re50Lift} and \ref{fig:Re50Drag}, the lift and drag coefficients are plotted, respectively. The drag was reduced by roughly 3.6\% for the QQR, while the vibration-inducing oscillations in the lift were reduced from a magnitude of 0.073 to 0.008 after 100 time units of actuation. Similarly to  the first case, the lift approaches zero as time goes to infinity.


\section{Conclusions and Future Work}

In conclusion, the pinball problem is a great model problem for testing flow control techniques. It exhibits strongly nonlinear behavior while remaining simple enough to simulate on modest computational resources. It also provides a useful framework for studying drag reduction and vortex-induced vibrations.

This paper introduces a control methodology that combines interpolatory model order reduction with the quadratic-quadratic regulator. The main result of this work is the finding that the QQR controller outperforms the linear controller for $Re_D=30$, and for $Re_D=50$, the QQR controller was able to stabilize the system while the linear controller failed. Moreover, it was found that the controller drives the lift fluctuations to zero while reducing the overall drag. Future work will involve a sensitivity analysis to investigate the impact of sensor placement. Additionally, further research could explore the influence of varying boundary conditions on flow dynamics and controller performance, as well as output feedback control. Understanding these aspects will help refine control strategies and enhance their applicability to more complex and realistic flow scenarios.

\clearpage
\bibliographystyle{plain}
\bibliography{pinball}

\end{document}